\documentclass[12pt,twoside,reqno,openany]{amsart}
\usepackage{courier}
\usepackage{etoolbox}
\makeatletter
\patchcmd{\@settitle}{\uppercasenonmath\@title}{}{}{}
\patchcmd{\@setauthors}{\MakeUppercase}{}{}{}
\patchcmd{\section}{\scshape}{}{}{}
\makeatother
\usepackage{amsbsy,amscd,amsfonts,amsmath,amssymb,amsthm,color,
fancybox,fancyhdr,footmisc,graphics,graphicx,ifthen,mathrsfs,
multicol,pdfpages,rotating,times,wasysym}

\usepackage[dvipsnames,svgnames,x11names]{xcolor}

\usepackage[all]{xy}
\usepackage[utf8]{inputenc}
\usepackage[T1]{fontenc}

\usepackage{mathtools}
\newtagform{EngelLie}[\scriptstyle]{$}{$}
\makeatletter\newcommand{\leqnomode}{\tagsleft@true}
\newcommand{\reqnomode}{\tagsleft@false}\makeatother

\newtheorem{Theorem}[equation]{Theorem}

\newtheorem{Proposition}[equation]{Proposition}

\newtheorem{Lemma}[equation]{Lemma}

\newtheorem{Corollary}[equation]{Corollary}

\theoremstyle{definition} 

\newtheorem{Hypothesis}[equation]{Hypothesis}
\newtheorem{Definition}[equation]{D\'efinition}

\newcommand{\FF}{\text{\sc f}}

\newcommand{\LL}{\text{\sc l}}

\newcommand{\NN}{\text{\sc n}}

\newcommand{\kaux}{{\text{\usefont{T1}{qcs}{m}{sl}k}}}

\newcommand{\Baux}{{\text{\usefont{T1}{qcs}{m}{sl}B}}}

\newcommand{\Iaux}{{\text{\usefont{T1}{qcs}{m}{sl}I}}}
\newcommand{\Jaux}{{\text{\usefont{T1}{qcs}{m}{sl}J}}}
\newcommand{\Kaux}{{\text{\usefont{T1}{qcs}{m}{sl}K}}}

\newcommand{\Paux}{{\text{\usefont{T1}{qcs}{m}{sl}P}}}
\newcommand{\Qaux}{{\text{\usefont{T1}{qcs}{m}{sl}Q}}}
\newcommand{\Raux}{{\text{\usefont{T1}{qcs}{m}{sl}R}}}

\newcommand{\Vaux}{{\text{\usefont{T1}{qcs}{m}{sl}V}}}
\newcommand{\Waux}{{\text{\usefont{T1}{qcs}{m}{sl}W}}}

\newcommand{\Zaux}{{\text{\usefont{T1}{qcs}{m}{sl}Z}}}

\definecolor{blue}{cmyk}{1.,1.,0.,0.63}
\definecolor{red}{cmyk}{0.,1.,1.,0.63}
\definecolor{green}{cmyk}{1.,0.,1.,0.63}
\definecolor{black}{cmyk}{1.,1.,1.,1.}

\makeatletter
\renewcommand{\@fnsymbol}[1]
{\ensuremath{\ifcase#1\or $*$\or $**$\or $***$\or $****$\or $*****$
\else\@ctrerr\fi}}
\makeatother

\numberwithin{equation}{section}

\newcommand{\style}[1]{\text{\footnotesize{\sf #1}}}

\renewcommand{\cosh}{\style{cosh}}

\renewcommand{\dim}{\style{dim}}

\renewcommand{\exp}{\style{exp}}

\renewcommand{\Im}{\style{Im}}

\renewcommand{\ker}{\style{ker}}

\newcommand{\Levi}{\style{Levi}}

\renewcommand{\lim}{\style{lim}}

\renewcommand{\mod}{\style{mod}}

\newcommand{\rank}{\style{rank}}

\renewcommand{\Re}{\style{Re}}

\renewcommand{\sinh}{\style{sinh}}

\renewcommand{\tanh}{\style{tanh}}

\newcommand{\isqrt}{i}

\newcommand{\LF}{\LL\FF}

\newcommand{\vf}{\vfill\end{document}}

\title[Equivalence problem of rigid biholomorphisms]{Rigid equivalences of $5$-dimensional $2$-nondegenerate\\ rigid real hypersurfaces $M^{5}\subset\mathbb{C}^{3}$ of constant Levi rank $1$\footnote{This work was supported in part by the Polish National Science Centre (NCN) via the grant number 2018/29/B/ST1/02583. The first author is supported by the Natural Science Foundation of China grant number 11688101.}}
\author{Wei Guo {\sc Foo}\qquad Jo\"el {\sc Merker}\qquad The-Anh {\sc Ta}}
\begin{document}
\begin{abstract}
We study the local equivalence problem for real-analytic ($\mathcal{C}^{\omega}$) hypersurfaces $M^{5}\subset\mathbb{C}^{3}$ which, in some holomorphic coordinates $(z_{1},z_{2},w)\in\mathbb{C}^{3}$ with $w=u+\isqrt v$, are \textit{rigid} in the sense that their graphing functions:
\[
u=F(z_{1},z_{2},\overline{z}_{1},\overline{z}_{2})
\]
are independent of $v$. Specifically, we study the group ${\sf Hol}_{\sf rigid}(M)$ of \textit{rigid} local biholomorphic transformations of the form: 
\[
\big(z_{1},z_{2},w\big)
\longmapsto 
\big(
f_{1}(z_{1},z_{2}),\ 
f_{2}(z_{1},z_{2}),\ 
aw+g(z_{1},z_{2})
\big),
\]
where $a\in \mathbb{R}\backslash \{0\}$ and $\frac{D(f_{1},f_{2})}{D(z_{1},z_{2})}\neq 0$, which preserve rigidity of hypersurfaces.

After performing a Cartan-type reduction to an appropriate $\{e\}$-structure, we find exactly \textit{two} primary invariants $\Iaux_{0}$ and $\Vaux_{0}$, which we express explicitly in terms of the $5$-jet of the graphing function $F$ of $M$. The identical vanishing $0\equiv \Iaux_{0}(J^{5}F)\equiv \Vaux_{0}(J^{5}F)$ then provides a necessary and sufficient condition for $M$ to be locally \textit{rigidly-biholomorphic} to the known model hypersurface:
\[
M_{\sf LC}:\qquad
u=
\frac{z_{1}\overline{z}_{1}+
\frac{1}{2}z_{1}^{2}\overline{z}_{2}
+
\frac{1}{2}\overline{z}_{1}^{2}z_{2}}{1-z_{2}\overline{z}_{2}}.
\]
We establish that $\dim\ {\sf Hol}_{\sf rigid}(M)\leqslant 7=\dim\ {\sf Hol}_{\sf rigid}(M_{\sf LC})$ always.

If one of these two primary invariants $\Iaux_{0}\not\equiv 0$ or $\Vaux_{0}\not\equiv 0$ does not vanish identically, then on either of the two Zariski-open sets $\{p\in M:\ \Iaux_{0}(p)\neq 0\}$ or $\{p\in M:\ \Vaux_{0}(p)\neq 0\}$, we show that this rigid equivalence problem between rigid hypersurfaces reduces to an equivalence problem for a certain $5$-dimensional $\{e\}$-structure on $M$, that is, we get an invariant absolute parallelism on $M^{5}$. Hence $\dim\ {\sf Hol}_{\sf rigid}(M)$ drops from $7$ to $5$, illustrating the \textit{gap phenomenon}.
\end{abstract}
\hspace{-1cm}
\maketitle
\tableofcontents

\section{Introduction}

In appropriate affine coordinates $(z_{1},z_{2},w)\in\mathbb{C}^{3}$ with $w=u+\isqrt v$, a real-analytic ($\mathcal{C}^{\omega}$) real hypersurface $M^{5}\subset\mathbb{C}^{3}$ may locally be represented as the graph of a $\mathcal{C}^{\omega}$ function $F$ over the $5$-dimensional real hyperplane $\mathbb{C}_{z_{1}}\times \mathbb{C}_{z_{2}}\times\mathbb{R}$. When $F$ is independent of $v$:
\[
M:\qquad u=F(z_{1},z_{2},\overline{z}_{1},\overline{z}_{2}),
\]
the hypersurface is called \textit{rigid}.

Its fundamental CR-bundle:
\[
T^{1,0}M\ :=\ 
\big(\mathbb{C}\otimes_{\mathbb{R}}TM)\cap T^{1,0}\mathbb{C}^{3}
\]
is of complex rank $2={\sf CRdim}M$, as well as its conjugate $T^{0,1}M=\overline{T^{1,0}M}$. 

Relevant foundational material for CR geometry focused on the local biholomorphic equivalence problem of $\mathcal{C}^{\omega}$ CR submanifolds $M\subset\mathbb{C}^{\NN}$ has be set up in the memoir \cite{MPS-2013}, to which readers will be referred for details.

The \textit{Levi forms} at various points $p\in M$ are maps measuring Lie bracket non-involutivity \cite[p. 45]{MPS-2013}:
\begin{equation*}
\begin{aligned}
T_{p}^{1,0}M\ \times\  T_{p}^{1,0}M 
&\longrightarrow 
\mathbb{C}\ \otimes_{\mathbb{R}}\ T_{p}M
\qquad
\mod\big(T_{p}^{1,0}M\oplus T_{p}^{0,1}M\big),\\
\big(\mathcal{M}_{p},\ \mathcal{N}_{p}\big)
&\longmapsto 
\isqrt \big[\mathcal{M},\ \overline{\mathcal{N}}\big]\big|_{p}\qquad
\mod\big(T_{p}^{1,0}M\oplus T_{p}^{0,1}M\big),
\end{aligned}
\end{equation*}
where $\mathcal{M}$ and $\mathcal{N}$ are any two local sections of $T^{1,0}M$ defined near $p$ which extend $\mathcal{M}_{p}=\mathcal{M}\big|_{p}$ and $\mathcal{N}_{p}=\mathcal{N}\big|_{p}$, the result being independent of extensions. 

Levi forms are known to be biholomorphically invariant. In terms of two natural intrinsic generators for $T^{1,0}M$:
\[
\mathcal{L}_{1}\ :=\ 
\frac{\partial}{\partial z_{1}}
-
\isqrt F_{z_{1}}\frac{\partial}{\partial v}
\qquad
\text{and}\qquad
\mathcal{L}_{2}\ :=\ 
\frac{\partial}{\partial z_{2}}
-\isqrt F_{z_{2}}\frac{\partial}{\partial v},
\]
the Levi forms at all points $p\in M$ identify with the matrix-valued map:
\[
\LF_{M}(p)
\ 
:=\ 
2\left(
\begin{matrix}
F_{z_{1}\overline{z}_{1}} & F_{z_{2}\overline{z}_{1}}\\
F_{z_{1}\overline{z}_{2}} & F_{z_{2}\overline{z}_{2}}
\end{matrix}
\right)(p).
\]
Throughout this article, we will make two main (invariant) assumptions. The first one is that the rank of $\LF_{M}(p)$ be constant equal to $1$ at every point $p\in M$.

Since $2=\rank\ T^{1,0}M$, this implies that there is a rank $1$ \textit{Levi kernel subbundle}:
\[
K^{1,0}M\ \subset\ T^{1,0}M,
\]
which is generated by the vector field:
\[
\mathcal{K}\ :=\ \kaux\mathcal{L}_{1}\ +\ \mathcal{L}_{2},
\]
incorporating the \textit{slant function}:
\[
\kaux\ :=\ 
-\frac{F_{z_{2}\overline{z}_{1}}}{F_{z_{1}\overline{z}_{1}}}.
\]
Indeed, a direct check convinces that both $[\mathcal{K},\ \overline{\mathcal{L}}_{1}]$ and $[\mathcal{K},\ \overline{\mathcal{L}}_{2}]$ vanish modulo $T^{1,0}M\oplus T^{0,1}M$. The known involutivity properties of the Levi kernerl subbundle $K^{1,0}M\subset T^{1,0}M$ together with its conjugate $K^{0,1}M\subset T^{0,1}M$ then read as (see \cite[pp. 72-73]{MPS-2013}):
\begin{equation*}
\begin{aligned}
\big[ K^{1,0}M,\ K^{1,0}M\big]\ &\subset\ K^{1,0}M,\\
\big[K^{0,1}M,\ K^{0,1}M\big]\ &\subset\ K^{0,1}M,\\
\big[K^{1,0}M,\ K^{0,1}M\big]\ &\subset\ K^{1,0}M\ \oplus\ K^{0,1}M.
\end{aligned}
\end{equation*}
Another fundamental function will also be needed in a while:
\[
\Paux\ :=\ \frac{F_{z_{1}z_{1}\overline{z}_{1}}}{F_{z_{1}\overline{z}_{1}}}.
\]

All this justifies the introduction of the so-called \textit{Freeman form} (\cite[p. 89]{MPS-2013}):
\begin{equation*}
\begin{aligned}
K_{p}^{1,0}M\ \times\ 
\big(T_{p}^{1,0}M\ \mod\ K_{p}^{1,0}M\big)
&\longrightarrow 
T_{p}^{1,0}M\ \oplus\ T_{p}^{0,1}M\quad 
\mod\ \big(K_{p}^{1,0}M\ \oplus\ T_{p}^{0,1}M\big),\\
\big(\mathcal{K}_{p},\ \mathcal{L}_{p}\big)
&\longmapsto 
\big[\mathcal{K},\ \overline{\mathcal{L}}\big]\big|_{p}
\quad\quad\quad\ \ 
\mod\ 
\big(K_{p}^{1,0}M\ \oplus\ T_{p}^{0,1}M\big),
\end{aligned}
\end{equation*}
where $\mathcal{K}$ and $\mathcal{L}$ are any two local sections of $K^{1,0}M$ and of $T^{1,0}M$ defined near $p$ which extend $\mathcal{K}_{p}=\mathcal{K}|_{p}$ and $\mathcal{L}_{p}=\mathcal{L}|_{p}$, the result being independent of extensions. In bases, these Freeman forms at various points $p\in M$ are simply maps $\mathbb{C}\times \mathbb{C}\longrightarrow\mathbb{C}$. They are known to be biholomorphically invariant. 

Our second main (invariant) assumption will be that the rank of the Freeman form be maximal equal to $1$ at every point $p\in M$. Such $M$ are called $2$-\textit{nondegenerate} at $p$.

A computation:
\begin{align*}
\big[\mathcal{K},\overline{\mathcal{L}}_{1}\big]=
\big[\kaux\mathcal{L}_{1}+\mathcal{L}_{2},\overline{\mathcal{L}}_{1}\big]
&=
-\overline{\mathcal{L}}_{1}(\kaux)\mathcal{L}_{1}
+
\underline{\kaux\ 
\big[\mathcal{L}_{1}, \overline{\mathcal{L}}_{1}\big]
+
\big[\mathcal{L}_{2},\overline{\mathcal{L}}_{1}\big]}_{\circ}\\
&=\ 
-\overline{\mathcal{L}}_{1}(\kaux)\mathcal{L}_{1}
\end{align*}
shows that 
\begin{center}
$M$ is $2$-nondegenerate at $p\in M$
\qquad
$\iff$
\qquad 
$\overline{\mathcal{L}}_{1}(\kaux)(p)\neq 0$.
\end{center}

Next, for a $\mathcal{C}^{\omega}$ hypersurface $M^{5}\subset\mathbb{C}^{3}$, define the Lie pseudogroup:
\[
{\sf Hol}_{\sf rigid}(M)
\ :=\ 
\bigg\{
h:\ M\longrightarrow M\ 
\text{local rigid biholomorphism}
\bigg\}.
\]
Its Lie algebra, obtained by differentiating $1$-parameter local groups of rigid biholomorphisms, is:
\begin{equation*}
\begin{aligned}
{\sf Lie}\big({\sf Hol}_{\sf rigid}(M)\big)
&=
\frak{hol}_{\sf rigid}(M)\\
&:= 
\bigg\{X=
A_{1}(z_{1},z_{2})\frac{\partial}{\partial z_{1}}
+
A_{2}(z_{1},z_{2})\frac{\partial}{\partial z_{2}}
+
(\alpha w+B(z_{1},z_{2}))
\frac{\partial}{\partial w}:
\\
&\hspace{1cm}
(X+\overline{X})|_{M}
\text{ is tangent to }M
\bigg\},
\end{aligned}
\end{equation*}
where $A_{1}$, $A_{2}$, $B$ are holomorphic functions of only $(z_{1},z_{2})$, and where $\alpha\in\mathbb{R}$. 

Our first result is the elementary
\begin{Proposition}
For the model hypersurface:
\[
M_{\sf LC}:\quad
u=
\frac{z_{1}\overline{z}_{1}+
\frac{1}{2}z_{1}^{2}\overline{z}_{2}
+
\frac{1}{2}\overline{z}_{1}^{2}z_{2}}{1-z_{2}\overline{z}_{2}},
\]
the Lie algebra $\frak{hol}_{\sf rigid}(M_{\sf LC})$ of infinitesimal biholomorphisms is $7$-dimensional, generated by:
\begin{equation*}
\begin{aligned}
X^{1} &= \isqrt\partial_{w},\\
X^{2} &= z_{1}\partial_{z_{1}} + 2w\partial_{w},\\
X^{3} &= \isqrt z_{1}\partial_{z_{1}} + 2\isqrt z_{2}\partial_{z_{2}},\\
X^{4} &= (z_{2}-1)\partial_{z_{1}} - 2z_{1}\partial_{w},\\
X^{5} &= (\isqrt+\isqrt z_{2})\partial_{z_{1}} - 2\isqrt z_{1}\partial_{w},\\
X^{6} &= z_{1}z_{2}\partial_{z_{1}} + (z_{2}^{2}-1)\partial_{z_{2}} - z_{1}^{2}\partial_{w},\\
X^{7} &= \isqrt z_{1}z_{2}\partial_{z_{1}} + (\isqrt z_{2}^{2}+\isqrt)\partial_{z_{2}} - \isqrt z_{1}^{2}\partial_{w}.
\end{aligned}
\end{equation*}
\end{Proposition}

Next, we conduct the Cartan process for rigid biholomorphic equivalences to this model $M_{\sf LC}$, reaching a representation of a Lie algebra isomorphic to the dual of the one generated by $X^{1},\dots,X^{7}$.

\begin{Theorem}
A basis for the Maurer-Cartan forms on the local Lie group ${\sf Hol}_{\sf rigid}(M_{\sf LC})$ is provided by $7$-differential $1$-forms:
\[
\{\rho,\ \kappa,\ \zeta,\ \overline{\kappa},\ \overline{\zeta},\ \alpha,\ \overline{\alpha}\},
\]
where $\overline{\rho}=\rho$ is real, which enjoys the $7$ structure equations with constant coefficients:
\begin{equation*}
\begin{aligned}
d\rho &= (\alpha+\overline{\alpha})\wedge\rho
+\isqrt\kappa\wedge\overline{\kappa},\\
d\kappa &= 
\alpha\wedge\kappa +\zeta\wedge\overline{\kappa},
\hspace{2cm}
d\overline{\kappa} = \overline{\alpha}\wedge\overline{\kappa}
+\overline{\zeta}\wedge\kappa,
\\
d\zeta &= (\alpha-\overline{\alpha})\wedge\zeta,
\hspace{2.4cm}
d\overline{\zeta} = (\overline{\alpha}-\alpha)\wedge\overline{\zeta},\\
d\alpha &= \zeta\wedge\overline{\zeta},
\hspace{3.5cm}
d\overline{\alpha} = \overline{\zeta}\wedge\zeta.
\end{aligned}
\qquad 
\end{equation*}
\end{Theorem}

This preliminary study of the model $M_{\sf LC}$ then constitutes our guiding map within the general problem. Recall that two fundamental functions expressed in terms of $F$ are:
\[
\kaux\ :=\ 
-\frac{F_{z_{2}\overline{z}_{1}}}{F_{z_{1}\overline{z}_{1}}}
\qquad
\text{and}
\qquad
\Paux\ :=\ 
\frac{F_{z_{1}z_{1}\overline{z}_{1}}}{F_{z_{1}\overline{z}_{1}}}.
\]

\begin{Theorem}\label{thm-main-equiv-2}
The equivalence problem under local rigid biholomorphisms of $\mathcal{C}^{\omega}$ rigid real hypersurfaces $\{u=F(z_{1},z_{2},\overline{z}_{1},\overline{z}_{2})\}$ in $\mathbb{C}^{3}$ whose Levi form has constant rank $1$ and which are everywhere $2$-nondegenerate reduces to classifying $\{e\}$-structures on the $7$-dimensional bundle $M^{5}\times \mathbb{C}$ equipped with coordinates $(z_{1},z_{2},\overline{z}_{1},\overline{z}_{2},v,{\sf c},\overline{\sf c})$ together with a coframe of $7$ differential $1$-forms:
\[
\{\rho,\ \kappa,\ \zeta,\ \overline{\kappa},\ \overline{\zeta},\ \alpha,\ \overline{\alpha}\},
\]
which satisfy invariant structure equations of the shape:
\begin{equation*}
\begin{aligned}
d\rho &= (\alpha+\overline{\alpha})\wedge\rho + \isqrt \kappa\wedge\overline{\kappa},\\
d\kappa &= \alpha\wedge\kappa + \zeta\wedge\overline{\kappa},\\
d\zeta &= (\alpha-\overline{\alpha})\wedge\zeta
+
\frac{1}{\sf c}\Iaux_{0}\ \kappa\wedge\zeta
+
\frac{1}{\overline{\sf c}\overline{\sf c}}
\Vaux_{0}\ \kappa\wedge\overline{\kappa},\\
d\alpha &= 
\zeta\wedge\overline{\zeta}
-
\frac{1}{\sf c}\Iaux_{0}\ \zeta\wedge\overline{\kappa}
+
\frac{1}{{\sf c}\overline{\sf c}}\Qaux_{0}\ \kappa\wedge\overline{\kappa}
+
\frac{1}{\overline{\sf c}}\overline{\Iaux}_{0}\ \overline{\zeta}\wedge\kappa,
\end{aligned}
\end{equation*}
conjugate equations for $d\overline{\kappa}$, $d\overline{\zeta}$, $d\overline{\alpha}$ being understood.
\end{Theorem}

Exactly two invariants are primary:
\begin{equation*}
\begin{aligned}
\Iaux_{0} &:= 
-\frac{1}{3}
\frac{\mathcal{K}(\overline{\mathcal{L}}_{1}(\overline{\mathcal{L}}_{1}(\kaux)))}{\overline{\mathcal{L}}_{1}(\kaux)^{2}}
+
\frac{1}{3}
\frac{\mathcal{K}(\overline{\mathcal{L}}_{1}(\kaux))
\overline{\mathcal{L}}_{1}(\overline{\mathcal{L}}_{1}(\kaux))}{\overline{\mathcal{L}}_{1}(\kaux)^{3}}
\\
&\hspace{0.5cm}
+
\frac{2}{3}
\frac{\mathcal{L}_{1}(\mathcal{L}_{1}(\overline{\kaux}))}{\mathcal{L}_{1}(\overline{\kaux})}
+\frac{2}{3}
\frac{\mathcal{L}_{1}(\overline{\mathcal{L}_{1}}(\kaux))}{\overline{\mathcal{L}}_{1}(\kaux)}
,\\
\Vaux_{0} &:=
-\frac{1}{3}
\frac{\overline{\mathcal{L}}_{1}(\overline{\mathcal{L}}_{1}(\overline{\mathcal{L}}_{1}(\kaux)))}{\overline{\mathcal{L}}_{1}(\kaux)}
+
\frac{5}{9}
\bigg(
\frac{\overline{\mathcal{L}_{1}}(\overline{\mathcal{L}}_{1}(\kaux))}{\overline{\mathcal{L}}_{1}(\kaux)}
\bigg)^{2}-\\
&\hspace{0.5cm}
-\frac{1}{9}
\frac{\overline{\mathcal{L}}_{1}(\overline{\mathcal{L}}_{1}(\kaux))\overline{\Paux}}{\overline{\mathcal{L}}_{1}(\kaux)}
+
\frac{1}{3}\overline{\mathcal{L}}_{1}(\overline{\Paux})
-
\frac{1}{9}
\overline{\Paux}\overline{\Paux},
\end{aligned}
\end{equation*}
while one invariant, which is real valued (see equation \eqref{Q0isreal}), is secondary:
\[
\Qaux_{0}
:=
\frac{1}{2}
\overline{\mathcal{L}}_{1}(\Iaux_{0})
-
\frac{1}{3}
\bigg(
\Paux
-
\frac{\mathcal{L}_{1}(\mathcal{L}_{1}(\overline{\kaux}))}{\mathcal{L}_{1}(\overline{\kaux})}
\bigg)
\overline{\Iaux_0}
-
\frac{1}{6}
\bigg(
\overline{\Paux}
-
\frac{\overline{\mathcal{L}}_{1}(\overline{\mathcal{L}}_{1}(\kaux))}{\overline{\mathcal{L}}_{1}(\kaux)}\bigg)\Iaux_{0}
-
\frac{1}{2}
\frac{\mathcal{K}(\Vaux_{0})}{\overline{\mathcal{L}}_{1}(\kaux)}.
\]

It is elementary to verify that both $\Iaux_{0}$ and $\Vaux_{0}$ vanish identically for $M_{\sf LC}$. As is known in Cartan theory, the identical vanishing of all invariants provide constant coefficients Maurer-Cartan equations of a uniquely defined Lie group. Hence as a corollary, we obtain the 

\begin{Theorem}
A 2-nondegenerate $\mathcal{C}^{\omega}$ constant Levi rank $1$ local rigid hypersurface $M^{5}\subset\mathbb{C}^{3}$ is rigidly biholomorphic to the model $M_{\sf LC}$ if and only if 
\[
0\equiv 
\Iaux_{0}\equiv \Vaux_{0}.
\eqno\qed
\]
\end{Theorem}

Next, when either $\Iaux_{0}\not\equiv 0$ or $\Vaux_{0}\not\equiv 0$, we may restrict considerations to either of the Zariski-open subsets $\{p\in M:\ \Iaux_{0}(p)\neq 0\}$ or $\{p\in M:\ \Vaux_{0}(p)\neq 0\}$, we may pursue the Cartan process, and we obtain the 

\begin{Theorem}
Let $M^{5}\subset \mathbb{C}^{3}$ be a local rigid $2$-nondegenerate $\mathcal{C}^{\omega}$ constant Levi rank $1$ hypersurface. If either $\Iaux_{0}\neq 0$ or $\Vaux_{0}\neq 0$ everywhere on $M$, the local rigid-biholomorphic equivalence problem reduces to an invariant $5$-dimensional $\{e\}$-structure on $M$.
\end{Theorem}

In fact, once the last remaining group parameter ${\sf c}\in \mathbb{C}^{*}$ is seen to be normalizable from either:
\[
\frac{1}{\sf c}\Iaux_{0}=1
\qquad
\text{or}
\qquad
\frac{1}{\overline{\sf c}\overline{\sf c}}\Vaux_{0}=1,
\]
the proof is completed if one does not require to make explicit the $\{e\}$-structure on $M$. Because of the size of computations, we will not attempt to set up such an explicit $\{e\}$-structure.

Lastly, from general Cartan theory, we deduce the 

\begin{Corollary}
All such rigid $M^{5}\subset\mathbb{C}^{3}$ that are not rigidly-biholomorphic to the model $M_{\sf LC}$ satisfy
\[
\dim\ {\sf Hol}_{\sf rigid}(M)\leqslant 5.
\]
\end{Corollary}

In continuation with these results, a further task appears: to classify up to rigid biholomorphisms the `\textit{submaximal}' hypersurfaces with $\dim\ {\sf Hol}_{\sf rigid}(M)=5$ whose rigid biholomorphic group is locally transitive. Another question would be to classify under rigid biholomorphisms those rigid $M^{5}\subset\mathbb{C}^{3}$ that have identically vanishing Pocchiola invariants $0\equiv \Waux_{0}\equiv \Jaux_{0}$, hence which are equivalent to $M_{\sf LC}$, but under a general biholomorphism, not necessarily rigid. Upcoming publications will be devoted to advances in these directions.


\section{Recall on the geometry of CR real hypersurfaces}

 Let $(z_{1},z_{2},w)$ be holomorphic coordinates in $\mathbb{C}^{3}$ with $w=u+\isqrt v$, and let $M^{5}\subset\mathbb{C}^{3}$ be a real-analytic, real hypersurface passing through the origin. Assuming that the real hypersurface is smooth at the origin, and that the vector 
 \[
 \frac{\partial}{\partial u}\bigg|_{0}\notin T_{0}M
 \]
 does not lie in the vector subspace $T_{0}M\subset T_{0}\mathbb{C}^{3}$. The implicit function theorem therefore implies the existence of a real analytic (denoted by $\mathcal{C}^{\omega}$) graphing function such that $M^{5}$ is represented near the origin by 
  \[
 u=F(z_{1},z_{2},\bar{z}_{1},\bar{z}_{2},v).
 \]
 
 \begin{Definition}
 The smooth real hypersurface $M^{5}\subset\mathbb{C}^{3}$ is {\sl rigid} at the origin if $M^{5}$ may be represented by a graphing function $u=F(z_{1},z_{2},\bar{z}_{1},\bar{z}_{2})$, where the function $F$ is independent of $v$.
 \end{Definition}
 
 \begin{Hypothesis}
 In the rest of the article, we will assume that $M^{5}$ is rigid.
 \end{Hypothesis}
 
 The complexified tangent bundle $\mathbb{C}TM=TM\otimes_{\mathbb{R}}\mathbb{C}$ inherits from $\mathbb{C}T\mathbb{C}^{3}=T\mathbb{C}^{3}\otimes_{\mathbb{R}}\mathbb{C}$ two biholomorphically invariant complex vector bundles
 \[
 T^{1,0}M:=\mathbb{C}TM\cap T^{1,0}\mathbb{C}^{3},\qquad
 T^{0,1}M:= \mathbb{C}TM\cap T^{0,1}\mathbb{C}^{3}
 =
 \overline{T^{1,0}M}.
 \]
 The two vector fields
 \[
 \mathcal{L}_{1} :=
 \frac{\partial}{\partial z_{1}}+A^{1}\frac{\partial}{\partial v}
 \qquad
 \text{and}
 \qquad
 \mathcal{L}_{2} :=
 \frac{\partial}{\partial z_{2}}
 +
 A^{2}\frac{\partial}{\partial v},
 \]
 with 
 \[
 A^{1} :=-\isqrt F_{z_{1}}, \qquad \text{and}\qquad
 A^{2} :=-\isqrt F_{z_{2}},
 \]
 then form a $T^{1,0}M$ frame. The differential $1$-form
 \[
 \rho_{0}
 =
 dv+\isqrt F_{z_{1}}\ dz^{1}
 +\isqrt F_{z_{2}}\ dz^{2}
 -\isqrt \overline{F}_{\bar{z}_{1}}\ d\bar{z}^{1}
 -\isqrt \overline{F}_{\bar{z}_{2}}\ d\bar{z}^{2}
 \]
 has the kernel
 \[
 \ker\rho_{0}
 =
 \{\rho_{0}=0\}
 =
 T^{1,0}M\oplus T^{0,1}M.
 \]
 By a formula in Merker-Pocchiola-Sabzevari \cite{MPS-2013}, page 82, the Levi matrix is shown to be 
 \begin{equation}
 \begin{aligned}
 \Levi(M)
 &=
  \left(
 \begin{matrix}
 \rho_{0}\big(\isqrt [\mathcal{L}_{1},\overline{\mathcal{L}}_{1}]\big)
 &
  \rho_{0}\big(\isqrt [\mathcal{L}_{2},\overline{\mathcal{L}}_{1}]\big)\\
   \rho_{0}\big(\isqrt [\mathcal{L}_{1},\overline{\mathcal{L}}_{2}]\big)
   &
    \rho_{0}\big(\isqrt [\mathcal{L}_{2},\overline{\mathcal{L}}_{2}]\big)
    \end{matrix}
 \right)\\
 &=
 2\left(
 \begin{matrix}
 F_{z_{1}\overline{z}_{1}} & 
 F_{z_{2}\overline{z}_{1}}\\
 F_{z_{1}\overline{z}_{2}} & 
 F_{z_{2}\overline{z}_{2}}
 \end{matrix}
 \right),
 \end{aligned}
 \end{equation}
 which is not identically zero if $M$ is further assumed to be not Levi-flat. After a change of coordinates in the $(z_{1},z_{2})$ space, without loss of generality, 
 \[
 \rho_{0}\big(\isqrt[\mathcal{L}_{1},\overline{\mathcal{L}}_{2}]\big)=
 2F_{z_{1}\bar{z}_{1}}\neq 0
 \]
everywhere on $M$, and hence the vector field
 \[
 \mathcal{T}
 :=
 \isqrt [\mathcal{L}_{1},\overline{\mathcal{L}}_{1}]
 =
 2F_{z_{1},\bar{z}_{1}}\frac{\partial}{\partial v}
 :=
 \ell \frac{\partial}{\partial v}
 \]
vanishes nowhere on $M$.
 
 \subsection{The rank 1 hypothesis}
 
 We will also make a further

 \begin{Hypothesis}
 The smooth real-analytic (rigid) real-hypersurface $M^{5}$ is of constant Levi rank 1.
 \end{Hypothesis}
 
 With this hypothesis, the collection of $1$-dimensional kernels $K_{p}^{1,0}M$ of the Levi form at all points $p\in M$ spans a real-analytic sub-distribution of the $T^{1,0}M$ bundle
 \[
 K^{1,0}M\subset T^{1,0}M,
 \]
 satisfying the following inclusions
 \begin{equation*}
 \begin{aligned}
 [K^{1,0}M,\ K^{1,0}M] &\subset K^{1,0}M,\\
 [K^{0,1}M,\ K^{0,1}M] &\subset K^{0,1}M,\\
 [K^{1,0}M,\ K^{0,1}M] &\subset K^{1,0}M\oplus K^{0,1}M.
 \end{aligned}
 \end{equation*}
 To construct a generator $\mathcal{K}$ of the Levi kernel, introduce a {\sl slant function} $\kaux$ satisfying
 \[
 \left(
 \begin{matrix}
 F_{z_{1},\bar{z}_{1}} & 
 F_{z_{2},\bar{z}_{1}}\\
 F_{z_{1}\bar{z}_{2}} & 
 F_{z_{2}\bar{z}_{2}}
 \end{matrix}
 \right)
 \left(
\begin{matrix}
\kaux\\
1
\end{matrix}
 \right)
 =
 \left(
 \begin{matrix}
 0\\
 0
 \end{matrix}
 \right).
 \]
 The first equation then implies that 
 \[
 \kaux = 
 -\frac{F_{z_{2}\bar{z}_{1}}}{F_{z_{1}\bar{z}_{1}}}
 \]
 while the same $\kaux$ satisfies the second equation
 \[
 \kaux F_{z_{1}\bar{z}_{2}}+F_{z_{2}\bar{z}_{2}}=0
 \]
 trivially by using the vanishing determinant of the matrix. Then the Levi kernel sub-bundle $\mathcal{K}^{1,0}M\subset T^{1,0}M$ is of complex rank $1$ and is generated by the vector field 
 \[
 \mathcal{K}=\kaux\mathcal{L}_{1}+\mathcal{L}_{2}.
 \]
  The slant function enjoys the following property
 \begin{Proposition}[See Merker-Pocchiola-Sabzevari \cite{MPS-2013}]
 The smooth real-analytic (rigid) real hypersurface $M$ is $2$-nondegenerate in the sense of Freeman if and only if 
 \[
 \overline{\mathcal{L}}_{1}(\kaux)\neq 0
 \]
 everywhere on $M$.
 \end{Proposition}
 In the rigid case, a direct calculation shows that 
 \begin{equation}
 \begin{aligned}
 \mathcal{L}_{1}(\kaux) 
 &= 
- \frac{-F_{z_{1},\bar{z}_{1}}F_{z_{2}\bar{z}_{1}z_{1}}
 +
 F_{z_{2}\bar{z}_{1}}F_{z_{1}\bar{z}_{1}z_{1}}}{(F_{z_{1}\bar{z}_{1}})^{2}},\\
 \overline{\mathcal{L}}_{1}(\kaux) &=
 \frac{-F_{z_{1}\bar{z}_{1}}F_{z_{2}\bar{z}_{1}\bar{z}_{1}}+
 F_{z_{2}\bar{z}_{1}}F_{z_{1}\bar{z}_{1}\bar{z}_{1}}}{(F_{z_{1}\bar{z}_{1}})^{2}},\\
 \mathcal{T}(\kaux) &= 0.
 \end{aligned}
 \end{equation}
 Moreover, introduce the next fundamental function
 \[
 \Paux = \frac{\ell_{z_{1}}}{\ell}=
 \frac{F_{z_{1}\bar{z}_{1}z_{1}}}{F_{z_{1}\bar{z}_{1}}}.
 \]
 
 \begin{Lemma}[See Pocchiola \cite{Pocchiola-2013} or Foo-Merker \cite{Foo-Merker-2019}]
 The following $3$ functional identities hold on $M$:
 \begin{equation}
 \begin{aligned}
 \mathcal{K}(\bar{\kaux}) &\equiv 0,\\
 \mathcal{K}(\Paux) &\equiv -\Paux\mathcal{L}_{1}(\kaux)-\mathcal{L}_{1}(\mathcal{L}_{1}(\kaux)),\\
 \mathcal{K}(\overline{\Paux}) &\equiv -\Paux\overline{\mathcal{L}}_{1}(\kaux)-\overline{\mathcal{L}}_{1}(\mathcal{L}_{1}(\kaux)).
 \end{aligned}
 \end{equation}
 \end{Lemma}
 
 According to Pocchiola \cite{Pocchiola-2013} page 37, there are 10 Lie bracket identities
 \begin{equation}
 \begin{aligned}
 [\mathcal{T},\mathcal{L}_{1}] &\equiv -\Paux\mathcal{T},\\
 [\mathcal{T},\mathcal{K}] &\equiv \mathcal{L}_{1}(\kaux)\mathcal{T} + 0,\\
 [\mathcal{T},\overline{\mathcal{L}}_{1}] &\equiv -\overline{P}\mathcal{T},\\
 [\mathcal{T},\overline{\mathcal{K}}] &\equiv \overline{\mathcal{L}}_{1}(\overline{\kaux})\mathcal{T} + 0,\\
 [\mathcal{L}_{1},\mathcal{K}] &\equiv \mathcal{L}_{1}(\kaux)\mathcal{L}_{1}
 \end{aligned}
 \qquad
 \begin{aligned}
 [\mathcal{L}_{1},\overline{\mathcal{L}}_{1}] &\equiv \isqrt\mathcal{T},\\
 [\mathcal{L}_{1},\overline{\mathcal{K}}] &\equiv \mathcal{L}_{1}(\overline{\kaux})\overline{\mathcal{L}}_{1},\\
 [\mathcal{K},\overline{\mathcal{L}}_{1}] &\equiv -\overline{\mathcal{L}}_{1}(\kaux)\mathcal{L}_{1},\\
 [\mathcal{K},\overline{\mathcal{K}}] &\equiv 0,\\
 [\overline{\mathcal{L}}_{1},\overline{\mathcal{K}}] &\equiv \overline{\mathcal{L}}_{1}(\overline{\kaux})\overline{\mathcal{L}}_{1}.
 \end{aligned}
 \end{equation}
 where the "$+0$" is deliberately added to show the difference from the general case. The following 1-forms
 \begin{equation}\label{ini-1-form}
 \begin{aligned}
 \rho_{0} &= \frac{1}{\ell}\big(dv-A^{1}dz_{1}-A^{2}dz_{2}-\bar{A}^{1}d\bar{z}_{1}-\bar{A}^{2}d\bar{z}_{2}\big),\\
 \kappa_{0} &= dz_{1}-\kaux dz_{2},\\
 \zeta_{0} &= dz_{2},\\
 \bar{\kappa}_{0} &= d\bar{z}_{1}-\bar{\kaux}d\bar{z}_{2},\\
 \bar{\zeta}_{0} &= d\bar{z}_{2},
 \end{aligned}
 \end{equation}
 are, by a simple computation, dual to the corresponding vector fields $\mathcal{T}$, $\mathcal{L}_{1}$, $\mathcal{K}$, $\overline{\mathcal{L}}_{1}$, $\overline{\mathcal{K}}$. Using the Cartan-Lie formula which states that for any smooth vector fields $X$, $Y$ and any smooth $1$-form $\omega$, one has
 \[
 d\omega(X,Y) = X\omega(Y)-Y\omega(X)-\omega([X,Y]),
 \]
 the initial Darboux-Cartan structure is therefore obtained:
 \begin{equation}\label{ini-CD}
 \begin{aligned}
 d\rho_{0} &= \Paux\ \rho_{0}\wedge\kappa_{0}
 -\mathcal{L}_{1}(\kaux)\ \rho_{0}\wedge\zeta_{0}
 +\overline{\Paux}\ \rho_{0}\wedge\bar{\kappa}_{0}
 -\overline{\mathcal{L}}_{1}(\bar{\kaux})\ \rho_{0}\wedge\bar{\zeta}_{0}
 +\isqrt\kappa_{0}\wedge\bar{\kappa}_{0},\\
 d\kappa_{0} &= -\mathcal{L}_{1}(\kaux)\ \kappa_{0}\wedge\zeta_{0}+
 \overline{\mathcal{L}}_{1}(\kaux)\ \zeta_{0}\wedge\bar{\kappa}_{0},\\
 d\zeta_{0} &= 0.
 \end{aligned}
 \end{equation}
 Here, conjugate equations for $d\overline{\kappa}_{0}$ and for $d\overline{\zeta}_{0}$ are not written, as they can be immediately deduced.


\section{Initial $G$-structure for rigid equivalences \\ of rigid real hypersurfaces}

Our objective is to study absolute parallelism of rigid equivalences of rigid real hypersurfaces using Cartan method. We introduce the

\begin{Definition}\label{def-rigid-aut}
Two local rigid real hypersurfaces at the origin are {\sl rigidly equivalent} if there exists a biholomorphic map of the form 
\[
\varphi:\ (z_{1},z_{2},w)\mapsto 
(z_{1}', z_{2}', w'):= \big(f(z_{1},z_{2}), g(z_{1},z_{2}), aw+h(z_{1},z_{2})\big),
\]
for some $a\in\mathbb{R}^{\times}$, and local holomorphic functions $f$, $g$, $h$, transforming one hypersurface into the other.
\end{Definition}

To make sure that the definition makes sense, let $M'$ be another rigid real hypersurface in the target space of the form
\[
\frac{w'+\bar{w}'}{2}
-
F'(z_{1}',z_{2}',\bar{z}_{1}',\bar{z}_{2}')=0.
\]
Then the pullback by $\varphi$ of the defining function is
\[
0
=
a\frac{w+\bar{w}}{2}
+
\bigg(
\frac{1}{2}h(z_{1},z_{2})
+
\frac{1}{2}\bar{h}(\bar{z}_{1},\bar{z}_{2})
-
F'
\big(f(z_{1},z_{2}),
g(z_{1},z_{2}),
\bar{f}(\bar{z}_{1},\bar{z}_{2}),
\bar{g}(\bar{z}_{1}\bar{z}_{2})
\big)
\bigg)
\]
which is again a defining function of a rigid real hypersurface.

Since $\varphi$ is holomorphic, its differential $\varphi_{*}:\mathbb{C}T\mathbb{C}^{3}\rightarrow \mathbb{C}T\mathbb{C}^{3}$ stablises the holomorphic $(1,0)$ and the anti-holomorphic $(0,1)$ vector fields:
\begin{equation}
\begin{aligned}
\varphi_{*}T^{1,0}M &\subseteq T^{1,0}M,\\
\varphi_{*}T^{0,1}M &\subseteq T^{0,1}M.
\end{aligned}
\end{equation}
Furthermore, by the invariance of the Freeman forms, the pushforward maps $\varphi_{*}$ also respects the Levi kernel distributions
\[
\varphi_{*}K^{1,0}M\subset K^{1,0}M.
\]
Consequently, there exist functions ${\sf f}'$, ${\sf c}'$ and ${\sf e}'$ on $M'$ such that 
\begin{equation}
\begin{aligned}
\varphi_{*}(\mathcal{K}) &= {\sf f}'\mathcal{K}',\\
\varphi_{*}(\mathcal{L}_{1}) &= {\sf c}'\mathcal{L}_{1}' + {\sf e}'\mathcal{K}'.
\end{aligned}
\end{equation}
The difference with the articles of Pocchiola \cite{Pocchiola-2013}, Merker-Pocchiola \cite{Merker-Pocchiola-2018} and Foo-Merker \cite{Foo-Merker-2019} is that the rigid equivalence assumption made on the map $\varphi:M\rightarrow M'$ between two rigid real hypersurfaces greatly simplifies the initial $G$-structure, especially because $\varphi_{*}\mathcal{T}$ is a multiple of $\mathcal{T}'$ by a function that vanishes nowhere on $M'$. In fact, if $R(z_{1}', z_{2}', \bar{z}_{1}', \bar{z}_{2}',v')$ is any $\mathcal{C}^{\omega}$ function on $M'$, then by definition of the pushforward of a vector field, 
\begin{equation}
\begin{aligned}
(\varphi_{*}\mathcal{T})\bigg(R(z_{1}',z_{2}',\bar{z}_{1}',\bar{z}_{2}',v')\bigg)
&=
\mathcal{T}(R\circ \varphi)\\
&=
\ell\frac{\partial}{\partial v}\big(R(f(z_{1},z_{2}), g(z_{1},z_{2}), \overline{f}(\overline{z}_{1},\overline{z}_{2}),\overline{g}(\overline{z}_{1},\overline{z}_{2}),av+\Im(h(z_{1},z_{2})))\big)\\
&=
a\ell\frac{\partial R}{\partial v'}\circ \varphi\\
&=
a\frac{\ell}{\ell'\circ\varphi}\underbrace{\bigg(\ell'\circ\varphi\frac{\partial R}{\partial v'}\circ \varphi\bigg)}_{=(\mathcal{T}'R)\circ\varphi}\\
&=
a\frac{\ell}{\ell'\circ\varphi}(\mathcal{T}'R)\circ\varphi,
\end{aligned}
\end{equation}
whence 
\[
\varphi_{*}\mathcal{T} = \frac{a\ell}{\ell'}\mathcal{T}'.
\]
Hence, there exists a real-valued function ${\sf a}' $ nowhere vanishing on $M'$ such that 
\[
\varphi_{*}\mathcal{T}={\sf a}'\mathcal{T}'.
\]

In fact, this function is determined since 
\begin{equation}
\begin{aligned}
{\sf a}'\mathcal{T}'=\varphi_{*}\mathcal{T} &= 
\varphi_{*}\big(\isqrt [\mathcal{L}_{1},\overline{\mathcal{L}}_{1}]\big)\\
&=
\isqrt [\varphi_{*}\mathcal{L}_{1},\varphi_{*}\overline{\mathcal{L}}_{1}]\\
&=
{\sf c}'{\sf\bar{c}}'\isqrt [\mathcal{L}_{1}',\overline{\mathcal{L}}_{1}'].
\end{aligned}
\end{equation}
This implies that 
\[
{\sf a}'={\sf c}'{\sf \bar{c}}'.
\]

Summarising, we therefore have the following matrix
\begin{equation}
\varphi_{*}
\left(
\begin{matrix}
\mathcal{T}\\
\mathcal{L}_{1}\\
\mathcal{K}\\
\overline{\mathcal{L}}_{1}\\
\overline{\mathcal{K}}
\end{matrix}
\right)
=
\left(
\begin{matrix}
{\sf c}'\bar{\sf c}' & 0 & 0 & 0 & 0\\
0 & {\sf c}' & {\sf e}' & 0 & 0 \\
0 & 0 & {\sf f}' & 0 & 0 \\
0 & 0 & 0 & \bar{\sf c}' & \bar{\sf e}'\\
0 & 0 & 0 & 0 & \bar{\sf f}'
\end{matrix}
\right)
\left(
\begin{matrix}
\mathcal{T}'\\
\mathcal{L}_{1}'\\
\mathcal{K}'\\
\overline{\mathcal{L}_{1}'}\\
\overline{\mathcal{K}}'
\end{matrix}
\right).
\end{equation}
Taking transposition of the matrix, one obtains the pullback formula for the two coframes
\begin{equation}\label{ini-g-str}
\begin{aligned}
\varphi^{*}\left(
\begin{matrix}
\rho_{0}'\\
\kappa_{0}'\\
\zeta_{0}'\\
\bar{\kappa}_{0}'\\
\bar{\zeta}_{0}'
\end{matrix}
\right)
=
\left(
\begin{matrix}
{\sf c}'\bar{\sf c}' & 0 & 0 & 0 & 0\\
0 & {\sf c}' & 0 & 0 & 0\\
0 & {\sf e}' & {\sf f}' & 0 & 0\\
0 & 0 & 0 & \bar{\sf c}' & 0\\
0 & 0 & 0 &  \bar{\sf e}' & \bar{\sf f}'
\end{matrix}
\right)
\left(
\begin{matrix}
\rho_{0}\\
\kappa_{0}\\
\zeta_{0}\\
\bar{\kappa}_{0}\\
\bar{\zeta}_{0}
\end{matrix}
\right)
\end{aligned}
\end{equation}
In conclusion, for the rigid CR transformation between rigid CR real hypersurfaces, the initial $G$-structure is constituted by the following 5 by 5 matrices
\[
\left(
\begin{matrix}
{\sf c}\bar{\sf c} & 0 & 0 & 0 & 0\\
0 & {\sf c} & 0 & 0 & 0\\
0 & {\sf e} & {\sf f} & 0 & 0\\
0 & 0 & 0 & \bar{\sf c} & 0\\
0 & 0 & 0 &  \bar{\sf e}' & \bar{\sf f}
\end{matrix}
\right)
\]
with the free complex variables 
\[
{\sf c}, {\sf f}\in \mathbb{C}-\{0\},\qquad
\text{and}
\qquad
{\sf e}\in\mathbb{C}.
\]


\section{Cartan equivalence method for the model case}

Before starting the Cartan equivalence method for rigid equivalences of $\mathcal{C}^{\omega}$ smooth rigid real hypersurfaces, a study of the equivalence method for the model case is necessary to obtain a model $\{e\}$-structure, which will serve as a reference for the general case. Recall that the model case is the tube over the future light cone, denoted by Pocchiola's notation as ${\sf MLC}$, is locally defined by the the following rigid equation
\[
u=
\frac{z_{1}\bar{z}_{1}+\frac{1}{2}z_{1}^{2}\bar{z}_{2}+\frac{1}{2}\bar{z}_{1}^{2}z_{2}}{1-z_{2}\bar{z}_{2}}.
\]
The vector fields $\mathcal{L}_{1}$, $\mathcal{K}$, $\mathcal{\overline{L}}_{1}$, $\overline{\mathcal{K}}$, $\mathcal{T}$, which constitute a frame for the complexified tangent bundle of $M_{\sf LC}$, thus have the following expressions
\begin{equation}
\begin{aligned}
\mathcal{L}_{1} &= \frac{\partial}{\partial z_{1}}-
\isqrt \frac{\bar{z}_{1}+z_{1}\bar{z}_{2}}{1-z_{2}\bar{z}_{2}}\frac{\partial}{\partial v},\\
\mathcal{K} &=
-\frac{\bar{z}_{1}+z_{1}\bar{z}_{2}}{1-z_{2}\bar{z}_{2}}\frac{\partial}{\partial z_{1}}
+
\frac{\partial}{\partial z_{2}}
+
\frac{\isqrt}{2}
\frac{\bar{z}_{1}^{2}+2z_{1}\bar{z}_{1}\bar{z}_{2}+z_{1}^{2}\bar{z}_{2}^{2}}{(1-z_{2}\bar{z}_{2})^{2}}
\frac{\partial}{\partial v},\\
\mathcal{T} &= 
-\frac{2}{1-z_{2}\bar{z}_{2}}\frac{\partial}{\partial v},
\end{aligned}
\end{equation}
and the slant function is given by
\[
\kaux = -\frac{\bar{z}_{1}+z_{1}\bar{z}_{2}}{1-z_{2}\bar{z}_{2}}.
\]
The initial coframe according to Pocchiola (model case) \cite{Pocchiola-2014} has the form
\begin{equation}
\begin{aligned}
\rho_{0} &=
-\frac{\isqrt}{2}(\bar{z}_{1}+z_{1}\bar{z}_{2})\ dz_{1}
-\frac{\isqrt}{4}\frac{\bar{z}_{1}^{2}+2z_{1}\bar{z}_{1}\bar{z}_{2}+z_{1}^{2}\bar{z}_{2}^{2}}{1-z_{2}\bar{z}_{2}}\ dz_{2}\\
&\hspace{0.5cm}
+\frac{\isqrt}{2}(z_{1}+\bar{z}_{1}z_{2})\ d\bar{z}_{2}
+
\frac{\isqrt}{4}
\frac{z_{1}^{2}+2z_{1}\bar{z}_{1}z_{2}+\bar{z}_{1}^{2}z_{2}^{2}}{1-z_{2}\bar{z}_{2}}\ d\bar{z}_{2}
+
\frac{1}{2}(-1+z_{2}\bar{z}_{2})\ dv,\\
\kappa_{0} &= 
dz_{1}+\frac{\bar{z}_{1}+z_{1}\bar{z}_{2}}{1-z_{2}\bar{z}_{2}}\ dz_{2},\\
\zeta_{0} &= dz_{2},
\end{aligned}
\end{equation}
which then satisfy the following structure equations
\begin{equation}
\begin{aligned}
d\rho_{0} &= 
\frac{\bar{z}_{2}}{1-z_{2}\bar{z}_{2}}\ \rho_{0}\wedge\zeta_{0}
+
\frac{z_{2}}{1-z_{2}\bar{z}_{2}}\ \rho_{0}\wedge\bar{\zeta}_{0}
+
\isqrt \kappa_{0}\wedge\bar{\kappa}_{0},\\
d\kappa_{0} &= 
\frac{\bar{z}_{2}}{1-z_{2}\bar{z}_{2}}\ \kappa_{0}\wedge\zeta_{0}
-
\frac{1}{1-z_{2}\bar{z}_{2}}\ \zeta_{0}\wedge\bar{\kappa}_{0},\\
d\zeta_{0} &= 0.
\end{aligned}
\end{equation}

In the case of rigid biholomorphisms as previously explained, the transformation group, denoted by ${\sf g}$, acts on the coframe $(\rho_{0},\kappa_{0},\zeta_{0})$ by the matrix
\[
{\sf g}=\left(
\begin{matrix}
{\sf c}\bar{\sf c} & 0 & 0\\
0 & {\sf c} & 0\\
0 & {\sf e} & {\sf f}
\end{matrix}
\right)
\]
while ignoring the $T^{0,1*}M$ counterpart. Its inverse
\[
{\sf g}^{-1}
=
\left(
\begin{matrix}
\frac{1}{{\sf c}\bar{\sf c}} & 0 & 0\\
0 & \frac{1}{\sf c} & 0\\
0 & -\frac{\sf e}{{\sf c}{\sf f}} & \frac{1}{\sf f}
\end{matrix}
\right)
\]
provides the following Maurer-Cartan matrix of $1$-forms 
\[
d{\sf g}\cdot {\sf g}^{-1}
=
\left(
\begin{matrix}
\alpha+\bar{\alpha} & 0 & 0\\
0 & \alpha & 0\\
0 & \delta & \varepsilon
\end{matrix}
\right),
\]
where the 1-forms $\alpha$, $\delta$, and $\varepsilon$ take on the following expressions
\begin{equation}\label{MC-1}
\begin{aligned}
\alpha &= \frac{d{\sf c}}{\sf c},\\
\delta &= \frac{d{\sf e}}{\sf c} - \frac{\sf e}{\sf c}\frac{d{\sf f}}{\sf f},\\
\varepsilon &= \frac{d{\sf f}}{\sf f}.
\end{aligned}
\end{equation}
Hence after some computation
\begin{equation}
\begin{aligned}
d\rho &=
(\alpha+\bar{\alpha})\wedge\rho
-
\frac{{\sf e}\bar{z}_{2}}{{\sf cf}(1-z_{2}\bar{z}_{2})}\ \rho\wedge\kappa
+
\frac{\bar{z}_{2}}{{\sf f}(1-z_{2}\bar{z}_{2})}\ \rho\wedge\zeta\\
&\hspace{2.5cm}
-
\frac{\bar{\sf e}z_{2}}{\overline{\sf cf}(1-z_{2}\bar{z}_{2})}\ \rho\wedge\bar{\kappa}
+
\frac{z_{2}}{\bar{\sf f}(1-z_{2}\bar{z}_{2})}\ \rho\wedge\bar{\zeta}
+
\isqrt\kappa\wedge\bar{\kappa},\\
d\kappa
&=
\alpha\wedge\kappa+
\frac{\bar{z}_{2}}{{\sf f}(1-z_{2}\bar{z}_{2})}\ \kappa\wedge\zeta
+
\frac{{\sf ce}}{{\sf c}\bar{\sf c}{\sf f}(1-z_{2}\bar{z}_{2})}\ \kappa\wedge\bar{\kappa}
-
\frac{{\sf c}}{\bar{\sf c}{\sf f}(1-z_{2}\bar{z}_{2})}\ \zeta\wedge\bar{\kappa},\\
d\zeta &=
\delta\wedge\kappa
+\varepsilon\wedge\zeta
+
\frac{{\sf e}\bar{z}_{2}}{{\sf cf}(1-z_{2}\bar{z}_{2})}\ \kappa\wedge\zeta
+
\frac{{\sf e}^{2}}{{\sf c}\bar{\sf c}{\sf f}(1-z_{2}\bar{z}_{2})}\ \kappa\wedge\bar{\kappa}\\
&\hspace{3cm}
+
\frac{\sf e}{\bar{\sf e}{\sf f}(1-z_{2}\bar{z}_{2})}\ \zeta\wedge\bar{\kappa}.
\end{aligned}
\end{equation}
In the rest of the article, we will adopt the following order for the coefficients appearing in front of the $2$-forms: 
\begin{equation}\label{ord-coeff}
\aligned
\underset{{\text{\bf 1}}}{
\rho_0\wedge\kappa_0}
\ \ \ \ \ \ \ \ \ \
\underset{{\text{\bf 2}}}{
\rho_0\wedge\zeta_0}
\ \ \ \ \ \ \ \ \ \
\underset{{\text{\bf 3}}}{
\rho_0\wedge\overline{\kappa}_0}
\ \ \ \ \ \ \ \ \ \
\underset{{\text{\bf 4}}}{
\rho_0\wedge\overline{\zeta}_0}
\\
\underset{{\text{\bf 5}}}{
\kappa_0\wedge\zeta_0}
\ \ \ \ \ \ \ \ \ \
\underset{{\text{\bf 6}}}{
\kappa_0\wedge\overline{\kappa}_0}
\ \ \ \ \ \ \ \ \ \
\underset{{\text{\bf 7}}}{
\kappa_0\wedge\overline{\zeta}_0}
\\
\underset{{\text{\bf 8}}}{
\zeta_0\wedge\overline{\kappa}_0}
\ \ \ \ \ \ \ \ \ \
\underset{{\text{\bf 9}}}{
\zeta_0\wedge\overline{\zeta}_0}
\\
\underset{{\text{\bf 10}}}{
\overline{\kappa}_0\wedge\overline{\zeta}_0}.
\endaligned
\end{equation}
Therefore, the $2$-forms may be abbreviated as
\begin{equation}
\begin{aligned}
d\rho &= 
(\alpha+\bar{\alpha})\wedge\rho
+
R_{1}\ \rho\wedge\kappa
+
R_{2}\ \rho\wedge\zeta
+
R_{3}\ \rho\wedge\bar{\kappa}
+
R_{4}\ \rho\wedge\bar{\zeta}\\
&\hspace{2.6cm}
+\isqrt \kappa\wedge\bar{\kappa},\\
d\kappa &=
\alpha\wedge\kappa + K_{5}\ \kappa\wedge\zeta + K_{6}\ \kappa\wedge\bar{\kappa} + K_{7}\ \zeta\wedge\bar{\kappa},\\
d\zeta &= \delta\wedge\kappa + \varepsilon\wedge\zeta
+Z_{5}\ \kappa\wedge\zeta+  Z_{6}\ \kappa\wedge\bar{\kappa}+Z_{8}\ \zeta\wedge\bar{\kappa}.
\end{aligned}
\end{equation}
Observe that $R_{3}=\overline{R}_{1}$ and $R_{4}=\overline{R}_{2}$. We will then proceed with the absorption, which can be done by replacing $\alpha$, $\delta$ and $\varepsilon$ with the new Maurer-Cartan $1$-forms
\begin{equation}
\begin{aligned}
\alpha &= 
\hat{\alpha}
-x_{\rho}\rho-x_{\kappa}\kappa-x_{\zeta}\zeta-x_{\bar{\kappa}}\bar{\kappa}-x_{\bar{\zeta}}\bar{\zeta},\\
\delta &= 
\hat{\delta} -y_{\rho}\rho-y_{\kappa}\kappa-y_{\zeta}\zeta-y_{\bar{\kappa}}\bar{\kappa}-y_{\bar{\zeta}}\bar{\zeta},\\
\varepsilon &= \hat{\varepsilon} -z_{\rho}\rho-z_{\kappa}\kappa-z_{\zeta}\zeta-z_{\bar{\kappa}}\bar{\kappa}-z_{\bar{\zeta}}\bar{\zeta}.
\end{aligned}
\end{equation}
for certain unknowns $x_{\bullet}$, $y_{\bullet}$ and $z_{\bullet}$. Therefore the $2$-forms may be re-written as 
\begin{equation}
\begin{aligned}
d\rho &= 
(\hat{\alpha}+\bar{\hat{\alpha}})\wedge\rho
+
(R_{1}+x_{\kappa}+\bar{x}_{\bar{\kappa}})\ \rho\wedge\kappa
+
(R_{2}+x_{\zeta}+\bar{x}_{\bar{\zeta}})\ \rho\wedge\zeta\\
&\hspace{3cm}
+(R_{3}+x_{\bar{\kappa}}+\overline{x_{\kappa}})\ \rho\wedge\bar{\kappa}
+
(R_{4}+x_{\bar{\zeta}}+\overline{x_{\zeta}})\ \rho\wedge\bar{\zeta}\\
&\hspace{3cm}
+\isqrt \kappa\wedge\bar{\kappa},\\
d\kappa &=
\hat{\alpha}\wedge\kappa +
(K_{5}+x_{\zeta})\ \kappa\wedge\zeta
+(K_{6}+x_{\bar{\kappa}})\ \kappa\wedge\bar{\kappa}
-
x_{\rho}\ \rho\wedge\kappa\\
&\hspace{3cm}
+
x_{\bar{\zeta}}\ \kappa\wedge\bar{\zeta}
+
K_{7}\ \zeta\wedge\bar{\kappa},\\
d\zeta &= \hat{\delta}\wedge\kappa+\hat{\varepsilon}\wedge\zeta
-y_{\rho}\ \rho\wedge\kappa
-z_{\rho}\ \rho\wedge\zeta
+(Z_{5}+y_{\zeta}-z_{\kappa})\ \kappa\wedge\zeta\\
&\hspace{1.5cm}
+
(Z_{6}+y_{\bar{\kappa}})\ \kappa\wedge\bar{\kappa}
+
(Z_{8}+z_{\bar{\kappa}})\ \zeta\wedge\bar{\kappa}
+
y_{\bar{\zeta}}\ \kappa\wedge\bar{\zeta}
+
z_{\bar{\zeta}}\ \zeta\wedge\bar{\zeta}.
\end{aligned}
\end{equation}
This therefore leads to the following set of equations
\begin{equation}
\begin{aligned}
x_{\kappa}+\overline{x_{\bar{\kappa}}} &= 
-\frac{{\sf e}\bar{z}_{2}}{{\sf cf}(1-z_{2}\bar{z}_{2})},\\
x_{\zeta}+\overline{x_{\bar{\zeta}}} &=
\frac{\bar{z}_{2}}{{\sf f}(1-z_{2}\bar{z}_{2})},\\
x_{\bar{\kappa}}+\overline{x_{\kappa}} &=
-\frac{{\sf e}z_{2}}{\bar{\sf cf}(1-z_{2}\bar{z}_{2})},\\
x_{\bar{\zeta}}+\overline{x_{\zeta}} &=
\frac{z_{2}}{\bar{\sf f}(1-z_{2}\bar{z}_{2})},\\
x_{\zeta} &= -\frac{\bar{z}_{2}}{{\sf f}(1-z_{2}\bar{z}_{2})},
\end{aligned}
\qquad
\begin{aligned}
x_{\bar{\kappa}} &= -\frac{\sf{ce}}{{\sf c}\bar{\sf c}{\sf f}(1-z_{2}\bar{z}_{2})},\\
x_{\rho} &= 0,\\
x_{\bar{\zeta}} &= 0,\\
y_{\rho} &= 0,\\
z_{\rho} &= 0,\\
\end{aligned}
\qquad
\begin{aligned}
y_{\zeta}-z_{\kappa} &= -\frac{{\sf e}\bar{z}_{2}}{{\sf cf}(1-z_{2}\bar{z}_{2})},\\
y_{\bar{\kappa}} &= -\frac{{\sf e}^{2}}{{\sf c}\bar{\sf c}{\sf f}(1-z_{2}\bar{z}_{2})},\\
z_{\bar{\kappa}} &= 
-\frac{{\sf e}}{\bar{\sf c}{\sf f}(1-z_{2}\bar{z}_{2})},\\
y_{\bar{\zeta}} &= 0,\\
z_{\bar{\zeta}} &= 0.
\end{aligned}
\end{equation}

These equations have solutions which result in the absorption of all the torsions except $K_{7}$, and hence 
\begin{equation}
\begin{aligned}
d\rho &= (\hat{\alpha}+\bar{\hat{\alpha}})\wedge\rho + \isqrt \kappa\wedge\bar{\kappa},\\
d\kappa &= \hat{\alpha}\wedge\kappa -
\frac{{\sf c}}{\bar{\sf c}{\sf f}(1-z_{2}\bar{z}_{2})}\ \zeta\wedge\bar{\kappa},\\
d\zeta &= \hat{\delta}\wedge\kappa+\hat{\varepsilon}\wedge\kappa.
\end{aligned}
\end{equation}
As in Pocchiola (model case) \cite{Pocchiola-2014}, the essential torsion 
\[
-\frac{{\sf c}}{\bar{\sf c}{\sf f}(1-z_{2}\bar{z}_{2})}
\]
may be normalised to $1$ by making the following choice
\[
{\sf f} = -\frac{{\sf c}}{\bar{\sf c}(1-z_{2}\bar{z}_{2})}.
\]

With this normalisation being made, we proceed with the second loop of the Cartan's equivalence method. The new transformation group then becomes
\begin{equation}
\begin{aligned}
\left(
\begin{matrix}
\rho\\
\kappa\\
\zeta
\end{matrix}
\right)
&=
\left(
\begin{matrix}
{\sf c}\bar{\sf c} & 0 & 0\\
0 & {\sf c} & 0\\
0 & {\sf e} & \frac{\sf c}{\bar{\sf c}}
\end{matrix}
\right)
\left(
\begin{matrix}
\rho_{0}\\
\kappa_{0}\\
-\frac{1}{1-z_{2}\bar{z}_{2}}\zeta_{0}
\end{matrix}
\right)\\
&:=
\left(
\begin{matrix}
{\sf c}\bar{\sf c} & 0 & 0\\
0 & {\sf c} & 0\\
0 & {\sf e} & \frac{\sf c}{\bar{\sf c}}
\end{matrix}
\right)
\left(
\begin{matrix}
\rho_{0}\\
\kappa_{0}\\
\hat{\zeta}_{0}
\end{matrix}
\right),
\end{aligned}
\end{equation}
with a change of the base coframe $(\rho_{0},\kappa_{0},\zeta_{0})\mapsto (\rho_{0},\kappa_{0},\hat{\zeta}_{0})$ via 
\[
\hat{\zeta}_{0} := -\frac{1}{1-z_{2}\overline{z}_{2}}\zeta_{0}.
\]
According to Pocchiola (model case) \cite{Pocchiola-2014}, the $2$-forms become 
\begin{equation}
\begin{aligned}
d\rho_{0} &= -\bar{z}_{2}\ \rho_{0}\wedge\hat{\zeta}_{0} - z_{2}\ \rho_{0}\wedge\bar{\hat{\zeta}}_{0}
+
\isqrt \kappa_{0}\wedge\bar{\kappa}_{0},\\
d\kappa_{0} &= -\bar{z}_{2}\ \kappa_{0}\wedge\hat{\zeta}_{0}+\hat{\zeta}_{0}\wedge\bar{\kappa}_{0},\\
d\hat{\zeta}_{0} &= z_{2}\ \hat{\zeta}_{0}\wedge\bar{\hat{\zeta}}_{0}.
\end{aligned}
\end{equation}
Moreover, one has the following Maurer-Cartan matrix of $1$-forms
\begin{equation}
\begin{aligned}
\left(
\begin{matrix}
\alpha+\bar{\alpha} & 0 & 0\\
0 & \alpha & 0\\
0 & \delta & \alpha-\bar{\alpha}
\end{matrix}
\right),
\end{aligned}
\end{equation}
where 
\[
\alpha = \frac{d{\sf c}}{\sf c},\qquad
\text{and}
\qquad
\delta = \frac{d{\sf e}}{\sf c}-\frac{\sf e}{\sf c}\left(
\frac{d{\sf c}}{\sf c}-\frac{d\bar{\sf c}}{\bar{\sf c}}\right).
\]
A computation by hand gives 
\begin{equation*}
\begin{aligned}
d\rho &= (\alpha+\bar{\alpha})\wedge\rho
+
\bar{z}_{2}\frac{{\sf e}\bar{\sf c}}{{\sf c}^{2}}\ \rho\wedge\kappa
-
\bar{z}_{2}\frac{\bar{\sf c}}{\sf c}\ \rho\wedge\zeta
+
z_{2}\frac{\bar{\sf e}{\sf c}}{\bar{\sf c}^{2}}\ \rho\wedge\bar{\kappa}
-
z_{2}\frac{{\sf c}}{\bar{\sf c}}\ \rho\wedge\bar{\zeta}
+
\isqrt\ \kappa\wedge\bar{\kappa},\\
&= 
(\alpha+\bar{\alpha})\wedge\rho
+
R_{1}\ \rho\wedge\kappa
+
R_{2}\ \rho\wedge\bar{\zeta}
+
\overline{R}_{1}\ \rho\wedge\bar{\kappa}
+
\overline{R}_{2}\ \rho\wedge\bar{\zeta}
+
\isqrt \kappa\wedge\bar{\kappa},\\
d\kappa &=
\alpha\wedge\kappa -\bar{z}_{2}\frac{\bar{\sf c}}{\sf c}\ \kappa\wedge\zeta
-\frac{\sf e}{\sf c}\ \kappa\wedge\bar{\kappa}
+
\zeta\wedge\bar{\kappa}\\
&=
\alpha\wedge\kappa+
K_{5}\ \kappa\wedge\zeta
+K_{6}\ \kappa\wedge\bar{\kappa}
+
\zeta\wedge\bar{\kappa},\\
d\zeta &= \delta\wedge\kappa
+(\alpha-\bar{\alpha})\wedge\zeta
-
\bar{z}_{2}\frac{{\sf e}\bar{\sf c}}{{\sf c}^{2}}\ \kappa\wedge\zeta
+
\left(
-\frac{{\sf e}^{2}}{{\sf c}^{2}}+z_{2}\frac{{\sf e}\bar{\sf e}}{{\bar{\sf c}}^{2}}\right)
\ \kappa\wedge\bar{\kappa}\\
&\hspace{3cm}
+
\left(
\frac{\sf e}{\sf c}-z_{2}
\frac{\bar{\sf e}{\sf c}}{\bar{\sf c}^{2}}
\right)\ \zeta\wedge\bar{\kappa}
-
z_{2}\frac{\sf e}{\bar{\sf c}}\ \kappa\wedge\bar{\zeta}
+
\frac{z_{2}{\sf c}}{\bar{\sf c}}\ \zeta\wedge\bar{\zeta}\\
&=
\delta\wedge\kappa+(\alpha-\bar{\alpha})\wedge\zeta
+ Z_{5}\ \kappa\wedge\zeta
+Z_{6}\ \kappa\wedge\bar{\kappa}
+Z_{8}\ \zeta\wedge\bar{\kappa}
+Z_{7}\ \kappa\wedge\bar{\zeta}
+Z_{9}\ \zeta\wedge\bar{\zeta}.
\end{aligned}
\end{equation*}
Then we proceed with the absorption by setting 
\begin{equation*}
\begin{aligned}
\alpha &= \hat{\alpha} -x_{\rho}\rho-x_{\kappa}\kappa-x_{\zeta}\zeta-x_{\bar{\kappa}}\bar{\kappa}-x_{\bar{\zeta}}\bar{\zeta},\\
\delta &= \hat{\delta}-y_{\rho}\rho-y_{\kappa}\kappa-y_{\zeta}\zeta-y_{\bar{\kappa}}\bar{\kappa}-y_{\bar{\zeta}}\bar{\zeta},
\end{aligned}
\end{equation*}
and we obtain
\begin{equation*}
\begin{aligned}
d\rho &= 
(\hat{\alpha}+\bar{\hat{\alpha}})\wedge\rho
+
(R_{1}+x_{\kappa}+\overline{x_{\bar{\kappa}}})\ \rho\wedge\kappa
+
(R_{2}+x_{\zeta}+\overline{x_{\bar{\zeta}}})\ \rho\wedge\zeta\\
&\hspace{1cm}
+
(\overline{R}_{1}+x_{\bar{\kappa}}+\overline{x_{\kappa}})\ \rho\wedge\bar{\kappa}
+
(\overline{R}_{2}+x_{\bar{\zeta}}+\overline{x_{\zeta}})\ \rho\wedge\bar{\zeta}
+
\isqrt\kappa\wedge\bar{\kappa},\\
d\kappa &= \hat{\alpha}\wedge\kappa
-x_{\rho}\ \rho\wedge\kappa
+(K_{5}+x_{\zeta})\ \kappa\wedge\zeta
+(K_{6}+x_{\bar{\kappa}})\ \kappa\wedge\bar{\kappa}
+x_{\bar{\zeta}}\ \kappa\wedge\bar{\zeta}
+\zeta\wedge\bar{\kappa},\\
d\zeta &= \hat{\delta}\wedge\kappa +
(\hat{\alpha}-\bar{\hat{\alpha}})\wedge\zeta
-y_{\rho}\ \rho\wedge\kappa
+(\overline{x_{\rho}}-x_{\rho})\ \rho\wedge\zeta
+
(Z_{5}-x_{\kappa}+y_{\zeta}+\overline{x_{\bar{\kappa}}})\ \kappa\wedge\zeta\\
&\hspace{1cm}
+(Z_{6}+y_{\bar{\kappa}})\ \kappa\wedge\bar{\kappa}
+Z_{7}\ \kappa\wedge\bar{\zeta}
+(Z_{8}-\overline{x_{\kappa}}+x_{\bar{\kappa}})\ \zeta\wedge\bar{\kappa}
+(Z_{9}-\overline{x_{\zeta}}+x_{\bar{\zeta}})\ \zeta\wedge\bar{\zeta}.
\end{aligned}
\end{equation*}
This leads to another set of absorption equations
\begin{equation}
\begin{aligned}
x_{\kappa}+\overline{x_{\bar{\kappa}}} &= -\bar{z}_{2}\frac{{\sf e}\bar{\sf c}}{{\sf c}^{2}},\\
x_{\zeta}+\overline{x_{\bar{\zeta}}} &= \bar{z}_{2}\frac{\bar{\sf c}}{\sf c},\\
x_{\bar{\kappa}}+\overline{x_{\kappa}} &= -z_{2}\frac{\bar{\sf e}{\sf c}}{\bar{\sf c}^{2}},\\
x_{\bar{\zeta}}+\overline{x_{\zeta}} &= z_{2}\frac{\sf c}{\bar{\sf c}},\\
x_{\rho} &= 0,\\
x_{\zeta} &= \bar{z}_{2}\frac{\bar{\sf c}}{\sf c},\\
x_{\bar{\kappa}} &= -\frac{\sf e}{\sf c},
\end{aligned}
\qquad
\begin{aligned}
x_{\bar{\zeta}} &= 0,\\
y_{\rho} &= 0,\\
\overline{x_{\rho}}-x_{\rho} &= 0,\\
-x_{\kappa}+\overline{x_{\bar{\kappa}}}+y_{\zeta}
&=
-\bar{z}_{2}\frac{{\sf e}\bar{\sf c}}{{\sf c}^{2}},\\
y_{\bar{\kappa}} &= -\frac{{\sf e}^{2}}{{\sf c}^{2}}
+z_{2}+\frac{{\sf e}\bar{\sf e}}{\bar{\sf c}^{2}},\\
-\overline{x_{\kappa}}+x_{\bar{\kappa}} &=
-\frac{\sf e}{\sf c}+z_{2}\frac{\bar{\sf e}{\sf c}}{\bar{\sf c}^{2}},\\
-\overline{x_{\zeta}}+x_{\bar{\zeta}} &= -z_{2}\frac{\sf c}{\bar{\sf c}}.
\end{aligned}
\end{equation}
The following equations
\begin{equation}
\begin{aligned}
x_{\bar{\kappa}}+\overline{x_{\kappa}} &= -z_{2}\frac{\bar{\sf e}{\sf c}}{\bar{\sf c}^{2}},\\
x_{\bar{\kappa}} &= -\frac{\sf e}{\sf c},\\
x_{\bar{\kappa}}-\overline{x_{\kappa}} &= -\frac{\sf e}{\sf c}+z_{2}\frac{\bar{\sf e}{\sf c}}{\bar{\sf c}^{2}}
\end{aligned}
\end{equation}
force us to conclude that ${\sf e}=0$, which is consistent with  Pocchiola (model case) \cite{Pocchiola-2014}, page 146, where he sets
\[
{\sf d}=-\isqrt \frac{{\sf e}^{2}\bar{\sf c}}{2{\sf c}}.
\]
In our case, ${\sf d}=0$ due to our rigidity assumption, and thus we are led to ${\sf e}=0$. 

This new normalisation gives rise to the new transformation group,
\[
\left(
\begin{matrix}
\rho\\
\kappa\\
\zeta
\end{matrix}
\right)
=
\left(
\begin{matrix}
{\sf c}\bar{\sf c} & 0 & 0\\
0 & {\sf c} & 0\\
0 & 0 & \frac{\sf c}{\bar{\sf c}}
\end{matrix}
\right)
\left(
\begin{matrix}
\rho_{0}\\
\kappa_{0}\\
\hat{\zeta}_{0}
\end{matrix}
\right)
\]
with the new Maurer-Cartan matrix
\[
d{\sf g}\cdot {\sf g}^{-1}
=
\left(
\begin{matrix}
\alpha+\bar{\alpha} & 0 & 0\\
0 & \alpha & 0\\
0 & 0 & \alpha-\bar{\alpha}
\end{matrix}
\right)
\]
and the following $2$-forms
\begin{equation}
\begin{aligned}
d\rho &= (\alpha+\bar{\alpha})\wedge\rho
-\bar{z}_{2}\frac{\bar{\sf c}}{\sf c}\ \rho\wedge\zeta
-z_{2}\frac{\sf c}{\bar{\sf c}}\ \rho\wedge\bar{\zeta}
+\isqrt \kappa\wedge\bar{\kappa},\\
d\kappa &= \alpha\wedge\kappa -
\bar{z}_{2}\frac{\bar{\sf c}}{\sf c}\ \kappa\wedge\zeta
+\zeta\wedge\bar{\kappa},\\
d\zeta &= (\alpha-\bar{\alpha})\wedge\zeta
+z_{2}\frac{\sf c}{\bar{\sf c}}\ \zeta\wedge\bar{\zeta}.
\end{aligned}
\end{equation}

We will proceed with the absorption process by setting 
\[
\alpha =
\hat{\alpha}
-x_{\rho}\rho
-x_{\kappa}\kappa
-x_{\zeta}\zeta
-x_{\bar{\kappa}}\bar{\kappa}
-x_{\bar{\zeta}}\bar{\zeta},
\]
which leads to
\begin{equation}
\begin{aligned}
d\rho &= (\hat{\alpha}+\bar{\hat{\alpha}})\wedge\rho
+(x_{\zeta}+\overline{x_{\bar{\zeta}}}-\bar{z}_{2}\frac{\bar{\sf c}}{\sf c})\ \rho\wedge\zeta
+(x_{\bar{\zeta}}+\overline{x_{\zeta}}-z_{2}\frac{\sf c}{\bar{\sf c}})\ \rho\wedge\bar{\zeta}\\
&\hspace{2cm}
+(x_{\kappa}+\overline{x_{\bar{\kappa}}})\ \rho\wedge\kappa
+(x_{\kappa}+\overline{x_{\kappa}})\ \rho\wedge\bar{\kappa}
+\isqrt\kappa\wedge\bar{\kappa},\\
d\kappa &=
\hat{\alpha}\wedge\kappa
-x_{\rho}\ \rho\wedge\kappa
+(x_{\zeta}-\bar{z}_{2}\frac{\bar{\sf c}}{\sf c})\ \kappa\wedge\zeta
+x_{\bar{\kappa}}\ \kappa\wedge\bar{\kappa}
+x_{\bar{\zeta}}\ \kappa\wedge\bar{\zeta}
+\zeta\wedge\bar{\kappa},\\
d\zeta &= (\hat{\alpha}-\bar{\hat{\alpha}})\wedge\zeta
+(-x_{\rho}+\overline{x_{\rho}})\ \rho\wedge\zeta
+(-x_{\kappa}+\overline{x_{\bar{\kappa}}})\ \kappa\wedge\zeta
+(\overline{x_{\kappa}}-x_{\bar{\kappa}})\ \zeta\wedge\bar{\kappa}\\
&\hspace{2cm}
+
(x_{\bar{\zeta}}-\overline{x_{\zeta}}+z_{2}\frac{\sf c}{\bar{\sf c}})\ \zeta\wedge\bar{\zeta}.
\end{aligned}
\end{equation}

To remove all the torsions, one has to solve for $x_{\bullet}$ the following system of linear equations
\begin{equation}
\begin{aligned}
x_{\zeta}+\overline{x_{\bar{\zeta}}} &= \bar{z}_{2}\frac{\bar{\sf c}}{\sf c},\\
x_{\bar{\zeta}}+\overline{x_{\zeta}} &= z_{2}\frac{\sf c}{\bar{\sf c}},\\
x_{\kappa}+\overline{x_{\bar{\kappa}}} &= 0,\\
x_{\rho} &= 0,\\
x_{\bar{\kappa}} &= 0,\\
x_{\bar{\zeta}} &= 0,\\
x_{\zeta} &= \bar{z}_{2}\frac{\bar{\sf c}}{\sf c},\\
-x_{\rho}+\overline{x_{\rho}} &= 0,\\
-x_{\kappa}+\overline{x_{\bar{\kappa}}} &= 0,\\
\overline{x_{\kappa}} -x_{\bar{\kappa}} &= 0,\\
x_{\bar{\zeta}}-\overline{x_{\zeta}}+z_{2}\frac{\sf c}{\bar{\sf c}} &= 0.
\end{aligned}
\end{equation}
This time the solution set is unambiguous:
\[
x_{\rho}=0,\qquad
x_{\kappa}=0,\qquad
x_{\zeta}=\bar{z}_{2}\frac{\bar{\sf c}}{\sf c},
\qquad
x_{\bar{\kappa}}=0,
\qquad
x_{\bar{\zeta}}=0,
\]
and since the degree of indeterminacy is zero, Cartan's test tells us that there is no need for prolongation. The absorption takes place and we get 
\begin{equation}
\begin{aligned}
d\rho &= (\hat{\alpha}+\bar{\hat{\alpha}})\wedge\rho+\isqrt \kappa\wedge\bar{\kappa},\\
d\kappa &= \hat{\alpha}\wedge\kappa +\zeta\wedge\bar{\kappa},\\
d\zeta &= (\hat{\alpha}-\bar{\hat{\alpha}})\wedge\zeta.
\end{aligned}
\end{equation}
The $\{e\}$-structure is then completed by the following
\begin{Proposition}
One has $d\hat{\alpha}=\zeta\wedge\bar{\zeta}$.
\end{Proposition}
\begin{proof}
Applying the Poincar\'{e} derivative on both sides of the three equations above, we get
\begin{eqnarray}
(d\hat{\alpha} + d\bar{\hat{\alpha}})\wedge\rho &=& 0,\label{model-prop-1-eq-1}\\
(d\hat{\alpha}-\zeta\wedge\bar{\zeta})\wedge\kappa &=& 0,\label{model-prop-1-eq-2}\\
(d\hat{\alpha}-d\bar{\hat{\alpha}})\wedge\zeta &=& 0.\label{model-prop-1-eq-3}
\end{eqnarray}
By applying complex conjugation to both sides of the third equation, one has an additional relation
\begin{eqnarray}
(d\hat{\alpha}-d\bar{\hat{\alpha}})\wedge\bar{\zeta} &=& 0.\label{model-prop-1-eq-4}
\end{eqnarray}
In the second equation \eqref{model-prop-1-eq-2}, Cartan's lemma provides a $1$-form $A$ so that 
\[
d\hat{\alpha} = \zeta\wedge\bar{\zeta} +A\wedge\kappa.
\]
Hence in  \eqref{model-prop-1-eq-1},  \eqref{model-prop-1-eq-3} and \eqref{model-prop-1-eq-4},
\begin{equation*}
\begin{aligned}
(d\hat{\alpha}+d\bar{\hat{\alpha}})\wedge\rho
&=
A\wedge\kappa\wedge\rho +
\bar{A}\wedge\bar{\kappa}\wedge \rho =0,\\
(d\hat{\alpha}-d\bar{\hat{\alpha}})\wedge\zeta
&=
A\wedge\kappa\wedge\zeta-\bar{A}\wedge\bar{\kappa}\wedge\zeta=0,\\
(d\hat{\alpha}-d\bar{\hat{\alpha}})\wedge\bar{\zeta}
&=
A\wedge\kappa\wedge\bar{\zeta}-\bar{A}\wedge\bar{\kappa}\wedge\bar{\zeta}=0.
\end{aligned}
\end{equation*}
Wedging with $\zeta$ on both sides of the first equation, and by $\rho$ on the second, we get
\begin{equation*}
\begin{aligned}
A\wedge\kappa\wedge\rho\wedge\zeta+\bar{A}\wedge\bar{\kappa}\wedge\rho\wedge\zeta=0,\\
A\wedge\kappa\wedge\rho\wedge\zeta-
\bar{A}\wedge\bar{\kappa}\wedge\rho\wedge\zeta=0.
\end{aligned}
\end{equation*}
Similarly, wedging with $\bar{\zeta}$ on both sides of the first equation, and by $\rho$ on the third, it comes
\begin{equation*}
\begin{aligned}
A\wedge\kappa\wedge\rho\wedge\bar{\zeta}+\bar{A}\wedge\bar{\kappa}\wedge\rho\wedge\bar{\zeta}=0,\\
A\wedge\kappa\wedge\rho\wedge\bar{\zeta}-
\bar{A}\wedge\bar{\kappa}\wedge\rho\wedge\bar{\zeta}=0.
\end{aligned}
\end{equation*}
Therefore,
\begin{equation}
\begin{aligned}
A\wedge\kappa\wedge\rho\wedge\zeta &= 0,\\
A\wedge\kappa\wedge\rho\wedge\bar{\zeta} &= 0.
\end{aligned}
\end{equation}
This implies the existence of functions $f$ and $g$ with
\[
A=f\rho+g\kappa.
\]
Hence
\[
d\hat{\alpha} = \zeta\wedge\bar{\zeta}+f\rho\wedge\kappa.
\]
Substituting this into \eqref{model-prop-1-eq-3}, 
\[
0 =
(\zeta\wedge\bar{\zeta}+f\rho\wedge\kappa
-
\bar{\zeta}\wedge\zeta
-
\bar{f}\rho\wedge\bar{\kappa})\wedge\zeta
=
f\rho\wedge\kappa\wedge\zeta
-\bar{f}\rho\wedge\bar{\kappa}\wedge\zeta,
\]
 we conclude by linear independence of these 3-forms that $f=0$. 
\end{proof}

\subsection{Summary}
For the model case, there exists a coframe $(\rho, \kappa, \zeta,\alpha, \bar{\kappa},\bar{\zeta},\bar{\alpha})$ satisfying the following structure equations
\begin{equation}\label{MLC-MC-structure-equations}
\begin{aligned}
d\rho &= (\alpha+\bar{\alpha})\wedge\rho + \isqrt \kappa\wedge\bar{\kappa},\\
d\kappa &= \alpha\wedge\kappa + \zeta\wedge\bar{\kappa},\\
d\zeta &= (\alpha-\bar{\alpha})\wedge\zeta,\\
d\alpha &= \zeta\wedge\bar{\zeta},
\end{aligned}
\end{equation}
along with the conjugates $d\rho$, $d\bar{\kappa}$, $d\bar{\zeta}$ and $d\bar{\alpha}$. Observe that $\alpha$ cannot be purely imaginary as seen  during the absorption of the final Cartan process. This therefore constitutes the Maurer-Cartan constant coefficients equations for the 7 dimensional complex Lie algebra of automorphisms of the model light cone $M_{\sf LC}$. We will confirm that this is $\frak{aut}_{\sf CR}(M_{\sf LC})$, arguing by means of  vector fields.


\section{Representation by vector fields}

By a result of Gaussier-Merker \cite{Gaussier-Merker-2003, Gaussier-Merker-Erratum}, it is known that the Lie algebra of infinitesimal CR automorphisms of the tube over future light cone $M_{\sf LC}$
is generated by the following 10 holomorphic vector fields 
\begin{equation}
\begin{aligned}
X^{1} &=
\isqrt \partial_{w},\\
X^{2} &=
z_{1}\partial_{z_{1}}+2w\partial_{w},\\
X^{3} &= 
\isqrt z_{1}\partial_{z_{1}}+2\isqrt z_{2}\partial_{z_{2}},\\
X^{4} &= (z_{2}-1)\partial_{z_{1}}-2z_{1}\partial_{w},\\
X^{5} &=
(\isqrt + \isqrt z_{2})\partial_{z_{1}}-2\isqrt z_{1}\partial_{w},\\
X^{6} &=
 z_{1}z_{2}\partial_{z_{1}}+(z_{2}^{2}-1)\partial_{z_{2}}-z_{1}^{2}\partial_{w},\\
 X^{7} &=
 \isqrt z_{1}z_{2}\partial_{z_{1}}+(\isqrt z_{2}^{2}+\isqrt)\partial_{z_{2}}-\isqrt z_{1}^{2}\partial_{w},\\
 X^{8} &=
 \isqrt wz_{1}\partial_{z_{1}}-\isqrt z_{1}^{2}\partial_{z_{2}}+\isqrt w^{2}\partial_{w},\\
 X^{9} &=
 (z_{1}^{2}-wz_{2}-w)\partial_{z_{1}}+(2z_{1}z_{2}+2z_{1})\partial_{z_{2}}+2wz_{1}\partial_{w},\\
 X^{10} &=
 (-\isqrt z_{1}^{2}+\isqrt wz_{2}-\isqrt w)\partial_{z_{1}}
 +(-2\isqrt z_{1}z_{2}+2\isqrt z_{1})\partial_{z_{2}}-2\isqrt wz_{1}\partial_{w}.
\end{aligned}
\end{equation}

It can be shown that for each $1\leqslant i\leqslant 10$, the vector field $X^{i}+\overline{X^{i}}$ is tangent to $M_{\sf LC}$. The commutator table of these 10 vector fields is as follows.

\begin{center}
  \begin{tabular}{ | c |c|c |c|c|c|c|c|c|c| c| }
    \hline
    & $X^{1}$ & $X^{2}$ & $X^{3}$ & $X^{4}$ & $X^{5}$ & $X^6$ & $X^7$ & $X^8$ & $X^9$ & $X^{10}$ \\ \hline
    $X^{1}$ & 0 & $2X^{1}$ &0 & 0 & 0 & 0 & 0 & $-X^{2}$ & $-X^{5}$ & $-X^{4}$\\ \hline
    $X^{2}$ &  & 0 & 0 & $-X^{4}$ & $-X^{5}$ & 0 & 0 & $2X^{8}$ & $X^{9}$ & $X^{10}$ \\ \hline
    $X^{3}$ &  &  & 0 & $X^{5}$ & $-X^{4}$ & $2X^{7}$ & $-2X^{6}$ & 0 & $-X^{10}$ & $X^{9}$\\ \hline
    $X^{4}$ &  &  &  & 0 & $4X^{1}$ & $-X^{4}$ & $-X^{5}$ & $X^{10}$ & $2X^{6}-2X^{2}$ & $-2X^{7}+2X^{3}$\\ \hline
    $X^{5}$ &  &  &  &  & 0 & $X^{5}$ & $-X^{4}$ & $X^{9}$ & $2X^{7}+2X^{3}$ & $2X^{6}+2X^{2}$\\ \hline
    $X^{6}$ &  &  &  &  &  & 0 & $-2X^{3}$ & 0 & $-X^{9}$ & $X^{10}$\\ \hline
    $X^{7}$ &  &  &  &  &  &  & 0 & 0 & $X^{10}$ & $X^{9}$\\ \hline
    $X^{8}$ &  &  &  &  &  &  &  & 0 & 0 & 0\\ \hline
    $X^{9}$ &  &  &  &  &  &  &  &  & 0 & $4X^{8}$\\ \hline
    $X^{10}$ &  &  &  &  &  &  &  &  &  & 0\\ \hline
  \end{tabular}
\end{center}

It is therefore clear from the table above that the vector fields $X^{1},\dots,X^{7}$ generate a Lie sub-algebra, which we will denote by $\frak{h}$. Next, we are going to find out which among these 10 vector fields have integral curves that define local rigid automorphisms of $\mathbb{C}^{3}$ (in the sense of Definition \ref{def-rigid-aut}). 

Recall that an integral curve of a  vector field $X$ on $\mathbb{C}^{3}$ is the map
\[
\gamma:\ \mathbb{R}\rightarrow \mathbb{C}^{3}
\]
satisfying the following differential equation with initial condition:
\begin{equation}
\begin{aligned}
\frac{d\gamma}{dt}\bigg|_{\gamma(t)} &= X|_{\gamma(t)},\\
\gamma(0) &= p.
\end{aligned}
\end{equation}
Usually such an integral curve at $p$ is denoted by 
\[
\exp(tX)(p) := \gamma(t)
\eqno
({\scriptstyle \gamma(0)=p}).
\]
Due to the following identity
\[
\exp(-tX)\,\exp(tX)(p)=p,
\]
an integral curve therefore defines an automorphism of $\mathbb{C}^{3}$ for each fixed $t$:
\begin{equation}
\begin{aligned}
\exp(tX):\ \mathbb{C}^{3} & \longrightarrow \mathbb{C}^{3}\\
p &\longmapsto \exp(tX)(p).
\end{aligned}
\end{equation}
For notational ease, we will let $p_{1}$, $p_{2}$ and $p_{3}$ denote the coordinates of 
\[
\gamma(0)=(\gamma_{1}(0),\gamma_{2}(0),\gamma_{3}(0))=(p_{1},p_{2},p_{3}).
\]

\subsection{Vector field $X^{1}$}
Integral curve: 
\[
(\gamma{1}(t),\gamma_{2}(t),\gamma_{3}(t))=
(p_{1},p_{2},p_{3}+it).
\]
Therefore for each fixed $t$, the holomorphic map
\[
(z_{1},z_{2},w)\mapsto (z_{1},z_{2},w+it)
\]
is rigid (we see $it$ as a constant holomorphic function).

\subsection{Vector field $X^{2}$}
Integral curve: 
\[
(\gamma_{1}(t),\gamma_{2}(t),\gamma_{3}(t))
=
(e^{t}p_{1},p_{2},e^{2t}p_{3}).
\]
Then for each fixed $t$,  the holomorphic map
\[
(z_{1},z_{2},w)\mapsto 
(e^{t}z_{1},z_{2},e^{2t}w)
\]
is rigid. 

\subsection{Vector field $X^{3}$}
Integral curve:
\[
(\gamma_{1}(t),\gamma_{2}(t),\gamma_{3}(t))
=
(e^{\isqrt t}p_{1},e^{2\isqrt t}p_{2},p_{3}).
\]
Therefore for each fixed $t$, the holomorphic map
\[
(z_{1},z_{2},w)
\mapsto
(e^{\isqrt t}z_{1},e^{2\isqrt t}z_{2},w)
\]
is rigid. 

\subsection{Vector field $X^4$}
Integral curve:
\[
(\gamma_{1}(t),\gamma_{2}(t),\gamma_{3}(t))
=
((p_{2}-1)t+p_{1},p_{2},-(p_{2}-1)t^{2}-2p_{1}t+p_{3}).
\]
For each fixed $t$, the holomorphic map
\[
(z_{1},z_{2},w)
\mapsto 
((z_{2}-1)t+z_{1},z_{2},w-((z_{2}-1)t^{2}+2z_{1}t))
\]
is rigid. 

\subsection{Vector field $X^{5}$}
Integral curve:
\[
(\gamma_{1}(t),\gamma_{2}(t),\gamma_{3}(t))
=
(p_{1}+\isqrt(p_{2}+1)t,p_{2},p_{3}-2\isqrt p_{1}t+(p_{2}+1)t^{2}).
\]
For each fixed $t$, the holomorphic map
\[
(z_{1},z_{2},w)
\mapsto 
(z_{1}+\isqrt (z_{2}+1)t, z_{2}, w-2\isqrt z_{1}t+(z_{2}+1)t^{2})
\]
is therefore rigid.

\subsection{Vector field $X^{6}$}
The integral curve $(\gamma_{1}(t),\gamma_{2}(t), \gamma_{3}(t))$ is given by the following equations
\begin{equation}
\begin{aligned}
\gamma_{1}(t) &= 
\frac{2p_{1}(1+p_{2})e^{t}}{(1+p_{2})(1+p_{2}+e^{2t}(1-p_{2}))},\\
\gamma_{2}(t) &=
\frac{(1+p_{2})-e^{2t}(1-p_{2})}{(1+p_{2})+e^{2t}(1-p_{2})},\\
\gamma_{3}(t) &= p_{3}+\frac{p_{1}^{2}}{1-p_{2}}
-\frac{2p_{1}^{2}}{1-p_{2}}\frac{1}{(1+p_{2})+(1-p_{2})e^{2t}}.
\end{aligned}
\end{equation}
For each fixed $t$, the holomorphic map
\begin{equation*}
\begin{aligned}
(z_{1},z_{2},w) &\mapsto
\bigg(
\frac{2z_{1}(1+z_{2})e^{t}}{(1+z_{2})(1+z_{2}+e^{2t}(1-z_{2}))},
\frac{(1+z_{2})-e^{2t}(1-z_{2})}{(1+z_{2})+e^{2t}(1-z_{2})},\\
&\hspace{3cm}
w+\frac{z_{1}^{2}}{1-z_{2}}-\frac{2z_{1}^{2}}{1-z_{2}}\frac{1}{(1+z_{2})+e^{2t}(1-z_{2})}
\bigg)
\end{aligned}
\end{equation*}
is therefore rigid. 

\subsection{Vector field $X^{7}$}
The integral curve $(\gamma_{1}(t),\gamma_{2}(t), \gamma_{3}(t))$ is given by
\begin{equation*}
\begin{aligned}
\gamma_{1}(t) &= \frac{\isqrt p_{1}}{p_{2}\sinh(t)+\isqrt \cosh(t)},\\
\gamma_{2}(t) &= 
\frac{-p_{2}-\isqrt\tanh(t)}{\isqrt + p_{2}\tanh(t)},\\
\gamma_{3}(t) &= p_{3}+\frac{p_{1}^{2}\sinh(t)}{p_{2}\sinh(t)+\isqrt \cosh(t)}.
\end{aligned}
\end{equation*}
Hence for each fixed $t$, the holomorphic map 
\begin{equation*}
\begin{aligned}
(z_{1},z_{2},w)
\mapsto
\bigg(
\frac{iz_{1}}{z_{2}\sinh(t)+\isqrt \cosh(t)},
\frac{-z_{2}-\isqrt\tanh(t)}{\isqrt+z_{2}\tanh(t)},
w+\frac{z_{1}^{2}\sinh(t)}{z_{2}\sinh(t)+\isqrt\cosh(t)}\bigg)
\end{aligned}
\end{equation*}
is rigid. 

One can deduce directly from the table that the Lie algebra $\frak{h}$ is neither semi-simple nor reductive. Indeed, the Killing form applied to the first vector field vanishes
\[
{\sf trace}({\sf ad}(X^{1}) {\sf ad}(X^{j}))=0,
\eqno
({\scriptstyle j=1,\dots,7})
\]
and hence $\frak{h}$ is not semi-simple by Cartan's criterion. Moreover, suppose by means of   {\sl reductio ad absurdum} that $\frak{h}$ is reductive, then it has a decomposition
\[
\frak{h}=\frak{s}\oplus \frak{z}(\frak{h}),
\]
where $\frak{s}$ is a semi-simple Lie sub-algebra and $\frak{z}(\frak{h})$ is the centre of $\frak{h}$. But it is clear from the table that $\frak{h}$ has no element in the centre except the zero vector field, and hence 
\[
\frak{h}=\frak{s}
\]
so that  $\frak{h}$ is semi-simple, a contradiction. 

We will now proceed to establish a link between the Maurer-Cartan coframe 
\begin{equation}\label{MC-coframe-model}
(\rho,\kappa,\zeta,\alpha,\bar{\kappa},\bar{\zeta},\bar{\alpha})
\end{equation}
 appearing in the structure equations in the previous sections, and the vector fields $X^{1}, \dots,\ X^{7}$. In fact, let 
\[
\partial_{\rho},\ \partial_{\kappa},\ \partial_{\zeta},\ \partial_{\alpha},\ \partial_{\bar{\kappa}},\ \partial_{\bar{\zeta}},\ \partial_{\bar{\alpha}}
\]
be the right-invariant vector fields that are respective duals to the $1$-forms in equation \eqref{MC-coframe-model}, and let $\frak{h}'$ be the Lie algebra generated by these vector fields. In what follows, the link will be established by seeking a Lie algebra isomorphism
\begin{equation}
\begin{aligned}
\tau: \frak{h}\longrightarrow \frak{h}'
\end{aligned}
\end{equation}
between $\frak{h}$ and $\frak{h}'$. 

We make the following recall which can be found in Olver \cite{Olver-1995}, page 257. Consider a set of $1$-forms $\theta=\{\theta^{1},\dots,\theta^{m}\}$ on a manifold $M$ producing the fundamental structure equations
\[
d\theta^{i}=\sum_{1\leqslant j<k\leqslant m}\ T^{i}_{jk}\ \theta^{j}\wedge\theta^{k}
\eqno
({\scriptstyle i=1,\dots,m}).
\]
If $\partial_{\theta^{i}}$ are the vector fields dual to $\theta^{i}$, one has the following commutation relations
\[
\big[\partial_{\theta^{j}},\partial_{\theta^{k}}\big]=
-\sum_{i=1}^{m}\ T_{jk}^{i}\ \partial_{\theta^{i}}
\eqno
({\scriptstyle 1\leqslant i<j\leqslant m}).
\]
Following this formula, and if we adopt the order of indices
\[\rho<\kappa<\zeta<\alpha<\bar{\kappa}<\bar{\zeta}<\bar{\alpha},\]
the Maurer-Cartan structure equations in equation \eqref{MLC-MC-structure-equations} therefore provide the following commutator table of the vector fields:
\newline
\begin{center}
  \begin{tabular}{ | c |c|c |c|c|c|c|c|}
    \hline
    & $\partial_{\rho}$ & $\partial_{\kappa}$ & $\partial_{\zeta}$ & $\partial_{\alpha}$ & $\partial_{\bar{\kappa}}$ & $\partial_{\bar{\zeta}}$ & $\partial_{\bar{\alpha}}$ \\ \hline
    $\partial_{\rho}$ & $0$ & $0$ & $0$ & $\partial_{\rho}$ & $0$ & $0$ & $\partial_{\bar{\rho}}$  \\ \hline
    $\partial_{\kappa}$ & $0$ & $0$ & $0$ & $\partial_{\kappa}$ & $-\isqrt\partial_{\rho}$ & $\partial_{\bar{\kappa}}$ & 0\\ \hline
    $\partial_{\zeta}$ & $0$ & $0$ & $0$ & $\partial_{\zeta}$ & $-\partial_{\kappa}$ & $-\partial_{\alpha}+\partial_{\bar{\alpha}}$ & $-\partial_{\zeta}$ \\ \hline
    $\partial_{\alpha}$ & $-\partial_{\rho}$ & $-\partial_{\kappa}$ & $-\partial_{\zeta}$ & $0$ & $0$ & $\partial_{\bar{\zeta}}$ & 0 \\ \hline
    $\partial_{\bar{\kappa}}$ & $0$ & $\isqrt\partial_{\rho}$ & $\partial_{\kappa}$ & $0$ & $0$ & $0$ & $\partial_{\bar{\kappa}}$ \\ \hline
    $\partial_{\bar{\zeta}}$ & $0$ & $-\partial_{\bar{\kappa}}$ & $-\partial_{\bar{\alpha}}+\partial_{\alpha}$ &$-\partial_{\bar{\zeta}}$ &  $0$  & $0$ & $\partial_{\bar{\zeta}}$ \\ \hline
    $\partial_{\bar{\alpha}}$ & $-\partial_{\rho}$ & $0$ &$\partial_{\zeta}$ & $0$ & $-\partial_{\bar{\kappa}}$ & $-\partial_{\bar{\zeta}}$ & 0\\ \hline
     \end{tabular}
\end{center}

Let $W^{1},\dots,\ W^{7}$ be the vector fields defined by 
\begin{equation}
\begin{aligned}
W^{1} &:=-\frac{\isqrt}{2}\partial_{\rho},\\
W^{2} &:=\partial_{\alpha}+\partial_{\bar{\alpha}},\\
W^{3} &:= \partial_{\zeta}-\partial_{\bar{\zeta}},
\end{aligned}
\qquad
\begin{aligned}
W^{4} &:= \partial_{\kappa}-\partial_{\bar{\kappa}},\\
W^{5} &:= \partial_{\kappa}+\partial_{\bar{\kappa}},\\
W^{6} &:= \partial_{\zeta}+\partial_{\bar{\zeta}},\\
W^{7} &:= -\partial_{\alpha}+\partial_{\bar{\alpha}}.
\end{aligned}
\end{equation}
Using the commutator table above, one has the following table of Lie brackets of various vector fields $W^{i}$:
\newline
\begin{center}
\begin{tabular}{|c|c|c|c|c|c|c|c|}
\hline
 & $W^{1}$ & $W^{2}$ & $W^{3}$ &$W^{4}$ & $W^{5}$ & $W^{6}$ & $W^{7}$\\ \hline
$W^{1}$ & $0$ & $2W^{1}$ & $0$ & $0$  & $0$ & $0$ & $0$ \\ \hline
$W^{2}$ &   & $0$ & $0$ & $-W^{4}$  & $-W^{5}$ & $0$ & $0$ \\ \hline
$W^{3}$ &  &  & $0$ & $W^{5}$  & $-W^{4}$ & $2W^{7}$ & $-2W^{6}$ \\ \hline
$W^{4}$ &  &  &  & $0$  & $4W^{1}$ & $-W^{4}$ & $-W^{5}$ \\ \hline
$W^{5}$ &  &  &  &  & $0$ & $W^{5}$ & $-W^{4}$ \\ \hline
$W^{6}$ &  &  &  &  & & $0$ & $-2W^{3}$ \\ \hline
$W^{7}$ &  &  & &  & & & $0$ \\ \hline
\end{tabular}
\end{center}
which is the same as the commutator table of the vector fields $X^{1},\dots,X^{7}$. Therefore the map which sends for each $i=1,\dots,7$:
\begin{equation}
\begin{aligned}
\tau:\ \mathfrak{h} & \longrightarrow \mathfrak{h}'\\
X^{i} & \longmapsto \tau(X^{i}):= W^{i}
\end{aligned}
\end{equation}
defines a Lie algebra isomorphism. The following theorem summarises what has been done so far for the rigid automorphisms of the model case:

\begin{Theorem}
The set of infinitesimal rigid CR-automorphisms of the tube over the future light cone 
\[
{\sf MLC}:\hfill
(\Re z_{1})^{2}-(\Re z_{2})^{2}-(\Re z_{3})^{2}=0
\eqno
\Re z_{1}>0,
\]
is a 7-dimensional Lie sub-algebra of the set of all of its infinitesimal CR-automorphisms. A basis for the Maurer-Cartan forms of the infinitesimal rigid CR-automorphisms is provided by the 7 differential $1$-forms $\rho$, $\kappa$, $\zeta$, $\alpha$, $\bar{\kappa}$, $\bar{\zeta}$, $\bar{\alpha}$ on ${\sf MLC}\times \mathbb{C}$ which satisfy the following Maurer-Cartan equations:
\begin{equation}
\begin{aligned}
d\rho &= (\alpha+\bar{\alpha})\wedge\rho+\isqrt \kappa\wedge\bar{\kappa},\\
d\kappa &= \alpha\wedge\kappa+\zeta\wedge\bar{\kappa},\\
d\zeta &= (\alpha-\bar{\alpha})\wedge\zeta,\\
d\alpha &= \zeta\wedge\bar{\zeta},\\
d\bar{\kappa} &= \bar{\alpha}\wedge\bar{\kappa} + \bar{\zeta}\wedge\kappa,\\
d\bar{\zeta} &= -(\alpha-\bar{\alpha})\wedge\bar{\zeta},\\
d\bar{\alpha} &= -\zeta\wedge\bar{\zeta}.
\end{aligned}
\end{equation}

Moreover, if $\{\partial_{\rho},\partial_{\kappa},\partial_{\zeta},\partial_{\alpha},\partial_{\bar{\kappa}},\partial_{\bar{\zeta}},\partial_{\bar{\alpha}}\}$ is a set of right-invariant vector fields that are dual to the respective coframe $1$-forms $\{\rho,\kappa,\zeta,\alpha,\bar{\kappa},\bar{\zeta},\bar{\alpha}\}$, then there is an isomorphism of Lie algebras between the Lie algebra $\frak{h}'$ generated by these vector fields, and the Lie algebra of infinitesimal rigid automorphisms of the tube over the future light cone. \qed
\end{Theorem}


\section{The general case}

The previous theorem shows that the Maurer-Cartan form that we have obtained, together with the structure equations, give a good setup for the equivalence problem. Recall from equations \eqref{ini-1-form} and \eqref{ini-CD} that the Darboux-Cartan structure equations are given by the $1$-forms $\{\rho_{0},\kappa_{0},\zeta_{0}\}$  with 
 \begin{equation}
 \begin{aligned}
 d\rho_{0} &= \Paux\ \rho_{0}\wedge\kappa_{0}
 -\mathcal{L}_{1}(\kaux)\ \rho_{0}\wedge\zeta_{0}
 +\overline{\Paux}\ \rho_{0}\wedge\bar{\kappa}_{0}
 -\overline{\mathcal{L}}_{1}(\bar{\kaux})\ \rho_{0}\wedge\bar{\zeta}_{0}
 +\isqrt\kappa_{0}\wedge\bar{\kappa}_{0},\\
 d\kappa_{0} &= -\mathcal{L}_{1}(\kaux)\ \kappa_{0}\wedge\zeta_{0}+
 \overline{\mathcal{L}}_{1}(\kaux)\ \zeta_{0}\wedge\bar{\kappa}_{0},\\
 d\zeta_{0} &= 0.
 \end{aligned}
 \end{equation}
 In equation \eqref{ini-g-str}, the group transformation of the $(1,0)$ coframe is determined by the matrix
 \[
 \omega
 =
 \left(
 \begin{matrix}
 \rho\\
 \kappa\\
 \zeta
 \end{matrix}
 \right)
 =
 \left(
 \begin{matrix}
 {\sf c}\bar{\sf c} & 0 & 0\\
 0 & {\sf c} & 0\\
 0 & {\sf e} & {\sf f}
 \end{matrix}
 \right)
 \left(
 \begin{matrix}
 \rho_{0}\\
 \kappa_{0}\\
 \zeta_{0}
 \end{matrix}
 \right)
 :=
 {\sf g}\omega_{0}.
 \]
 We will also continue to adopt the order of coefficients as stated in equation \eqref{ord-coeff}
 
 \section{Cartan process: first loop}
 Using the formula
 \[
 d\omega=(d{\sf g}){\sf g}^{-1}\omega+{\sf g}d\omega_{0},
 \]
 the Maurer-Cartan form is 
 \[
 (d{\sf g}){\sf g}^{-1}
 =
 \left(
 \begin{matrix}
 \alpha+\bar{\alpha} & 0 & 0\\
 0 & \alpha & 0\\
 0 & \delta & \varepsilon
 \end{matrix}
 \right),
 \]
 where $\alpha$, $\delta$ and $\varepsilon$ are given by those in equation \eqref{MC-1}. A direct computation shows that 
 \begin{equation}
 \begin{aligned}
 d\rho &= \alpha\wedge\rho + \bar{\alpha}\wedge\rho
 +
 \bigg(
 \frac{\Paux}{\sf c}+\frac{{\sf e}\mathcal{L}_{1}(\kaux)}{\sf cf}
 \bigg)\ \rho\wedge\kappa
 +
 \bigg(
 \frac{\bar{\Paux}}{\sf \bar{c}}
 +
 \frac{\bar{\sf e}\overline{\mathcal{L}}_{1}(\bar{\kaux})}{\bar{\sf c}\bar{\sf f}}
 \bigg)\ \rho\wedge\bar{\kappa}\\
 &\hspace{0.5cm}
 +\bigg(
 \frac{-\mathcal{L}_{1}(\kaux)}{\sf f}
 \bigg)\ \rho\wedge\zeta
 +
 \bigg(\frac{-\overline{\mathcal{L}}_{1}(\bar{\kaux})}{\bar{\sf f}}\bigg)\ 
 \rho\wedge\bar{\zeta}
 +
 \isqrt\kappa\wedge\bar{\kappa},\\
 d\kappa &= 
 \alpha\wedge\kappa + 
 \bigg(\frac{-\mathcal{L}_{1}(\kappa)}{\sf f}\bigg)\ \kappa\wedge\zeta
 +
 \bigg(-\frac{{\sf e}\overline{\mathcal{L}}_{1}(\kaux)}{\bar{\sf c}{\sf f}}\bigg)\ 
 \kappa\wedge\bar{\kappa}
 +
 \bigg(\frac{{\sf c}\overline{\mathcal{L}}_{1}(\kaux)}{\bar{\sf c}{\sf f}}\bigg)\ \zeta\wedge\bar{\kappa},\\
 d\zeta &= 
 \delta\wedge\kappa
 +\varepsilon\wedge\zeta
 +
 \bigg(\frac{-{\sf e}\mathcal{L}_{1}(\kaux)}{\sf cf}\bigg)\ \kappa\wedge\zeta
 +
 \bigg(\frac{-{\sf e}^{2}\overline{\mathcal{L}}_{1}(\kaux)}{{\sf c}\bar{\sf c}{\sf f}}\bigg)\ \kappa\wedge\bar{\kappa}\\
 &\hspace{0.5cm}+
 \bigg(\frac{{\sf e}\mathcal{L}_{1}(\kaux)}{\bar{\sf c}{\sf f}}\bigg)\ \zeta\wedge\bar{\kappa}.
 \end{aligned}
 \end{equation}
We proceed with the absorption by setting 
\begin{equation*}
\begin{aligned}
\alpha &= 
\hat{\alpha}
-x_{\rho}\rho-x_{\kappa}\kappa-x_{\zeta}\zeta-x_{\bar{\kappa}}\bar{\kappa}-x_{\bar{\zeta}}\bar{\zeta},\\
\delta &= 
\hat{\delta} -y_{\rho}\rho-y_{\kappa}\kappa-y_{\zeta}\zeta-y_{\bar{\kappa}}\bar{\kappa}-y_{\bar{\zeta}}\bar{\zeta},\\
\varepsilon &= \hat{\varepsilon} -z_{\rho}\rho-z_{\kappa}\kappa-z_{\zeta}\zeta-z_{\bar{\kappa}}\bar{\kappa}-z_{\bar{\zeta}}\bar{\zeta}.
\end{aligned}
\end{equation*}
Solving a system of linear equations to eliminate as many torsions as possible, one obtains
\begin{equation}
\begin{aligned}
d\rho &= (\hat{\alpha}+\overline{\hat{\alpha}})\wedge\rho + \isqrt \kappa\wedge\bar{\kappa},\\
d\kappa &= \hat{\alpha}\wedge\kappa + 
\frac{{\sf c}\overline{\mathcal{L}}_{1}(\kaux)}{\bar{\sf c}{\sf f}}\ \zeta\wedge\bar{\kappa},\\
d\zeta &= \hat{\delta}\wedge\kappa+\hat{\varepsilon}\wedge\zeta.
\end{aligned}
\end{equation}
Notice that the function 
\[
\frac{{\sf c}\overline{\mathcal{L}}_{1}(\kaux)}{\bar{\sf c}}
\]
is nowhere vanishing, and hence the torsion that appears in $d\kappa$ may be normalised to $1$ by setting
\[
{\sf f}=\frac{{\sf c}\overline{\mathcal{L}}_{1}(\kaux)}{\bar{\sf c}}.
\]

\section{Cartan process: second loop}
With this normalisation, we proceed with a change of the base coframe
\[
\hat{\zeta}_{0}:= \overline{\mathcal{L}}_{1}(\kappa)\zeta_{0},
\]
so that the new transformation group becomes 
\[
\left(
\begin{matrix}
\rho\\
\kappa\\
\zeta
\end{matrix}
\right)
=
\left(
\begin{matrix}
\sf{c}\bar{\sf c} & 0 & 0\\
0 & {\sf c} & 0\\
0 & {\sf e} & \frac{\sf c}{\bar{\sf c}}
\end{matrix}
\right)
\left(
\begin{matrix}
\rho_{0}\\
\kappa_{0}\\
\hat{\zeta}_{0}
\end{matrix}
\right).
\]
Observe that both functions vanish identically
\[
\mathcal{T}(\kaux)\equiv 0,\qquad 
\mathcal{T}(\mathcal{L}_{1}(\kaux))\equiv 0,
\]
since both $\kaux$ and $\mathcal{L}_{1}(\kaux)$ are independent of $v$. Using equation (5.5) of Foo-Merker \cite{Foo-Merker-2019}, the new Darboux-Cartan structure equations become 
\begin{equation}
\begin{aligned}
d\rho_{0} &= \Paux\ \rho\wedge\kappa_{0}
-\frac{\mathcal{L}_{1}(\kaux)}{\overline{\mathcal{L}}_{1}(\kaux)}\ \rho\wedge\hat{\zeta}_{0}
+\overline{\Paux}\ \rho_{0}\wedge\bar{\kappa}_{0}
-\frac{\overline{\mathcal{L}}_{1}(\bar{\kaux})}{\mathcal{L}_{1}(\bar{\kaux})}\ \rho_{0}\wedge\overline{\hat{\zeta}_{0}}+\isqrt\kappa_{0}\wedge\bar{\kappa}_{0},\\
d\kappa_{0} &= 
-\frac{\mathcal{L}_{1}(\kaux)}{\overline{\mathcal{L}}_{1}(\kaux)}\ \kappa_{0}\wedge\hat{\zeta}_{0}
+\hat{\zeta}_{0}\wedge\bar{\kappa}_{0},\\
d\hat{\zeta}_{0} &= 
\frac{\mathcal{L}_{1}(\overline{\mathcal{L}}_{1}(\kaux))}{\overline{\mathcal{L}}_{1}(\kaux)}\ \kappa_{0}\wedge\hat{\zeta}_{0}-
\frac{\overline{\mathcal{L}}_{1}(\overline{\mathcal{L}}_{1}(\kaux))}{\overline{\mathcal{L}}_{1}(\kaux)}\ 
\hat{\zeta}_{0}\wedge\overline{\kappa}_{0}
+\frac{\overline{\mathcal{L}}_{1}(\bar{\kaux})}{\mathcal{L}_{1}(\bar{\kaux})}\ \hat{\zeta}_{0}\wedge\overline{\hat{\zeta}_{0}}.
\end{aligned}
\end{equation}
Moreover, one has the following Maurer-Cartan matrix
\[
(d{\sf g}){\sf g}^{-1}
=
\left(
\begin{matrix}
\alpha+\bar{\alpha} & 0 & 0\\
0 & \alpha & 0\\
0 & \delta & \alpha-\bar{\alpha}
\end{matrix}
\right),
\]
with the 1-forms
\[
\alpha=\frac{d{\sf c}}{\sf c},\qquad 
\delta=\frac{d{\sf e}}{\sf c}-\frac{\sf e}{\sf c}
\bigg(\frac{d{\sf c}}{\sf c}-\frac{d\bar{\sf c}}{\bar{\sf c}}\bigg).
\]
One obtains therefore
\begin{equation*}
\begin{aligned}
d\rho &= (\alpha+\bar{\alpha})\wedge\rho
+
\bigg(\frac{\Paux}{\sf c}+
\frac{\mathcal{L}_{1}(\kaux)}{\overline{\mathcal{L}}_{1}(\kaux)}
\frac{{\sf e}\bar{\sf c}}{{\sf c}^{2}}\bigg)\ \rho\wedge\kappa
+
\bigg(
-\frac{\mathcal{L}_{1}(\kaux)}{\overline{\mathcal{L}}_{1}(\kaux)}\frac{\bar{\sf c}}{\sf c}\bigg)\ 
\rho\wedge\zeta\\
&\hspace{3cm}
+\bigg(\frac{\bar{\Paux}}{\bar{\sf c}}
+\frac{\overline{\mathcal{L}}_{1}(\bar{\kaux})}{\mathcal{L}_{1}(\bar{\kaux})}
\frac{\bar{\sf e}{\sf c}}{\bar{\sf c}^{2}}\bigg)\ \rho\wedge\bar{\kappa}
+
\bigg(-\frac{\overline{\mathcal{L}}_{1}(\bar{\kaux})}{\mathcal{L}_{1}(\bar{\kaux})}\frac{\sf c}{\bar{\sf c}}\bigg)\ \rho\wedge\bar{\zeta}+\isqrt\kappa\wedge\bar{\kappa},\\
d\kappa &=
\alpha\wedge\kappa
+\bigg(-\frac{\mathcal{L}_{1}(\kaux)}{\overline{\mathcal{L}}_{1}(\kaux)}\frac{\bar{\sf c}}{\sf c}\bigg)\ 
\kappa\wedge\zeta-
\frac{\sf e}{\sf c}\ \kappa\wedge\bar{\kappa}
+
\zeta\wedge\bar{\kappa},\\
d\zeta &= 
\delta\wedge\kappa +
(\alpha-\bar{\alpha})\wedge\zeta
+
\bigg(-\frac{\mathcal{L}_{1}(\kaux)}{\overline{\mathcal{L}}_{1}(\kaux)}\frac{{\sf e}\bar{\sf c}}{{\sf c}^{2}}
+\frac{\mathcal{L}_{1}(\overline{\mathcal{L}}_{1}(\kaux))}{\overline{\mathcal{L}}_{1}(\kaux)}\frac{1}{\sf c}\bigg)\ \kappa\wedge\zeta\\
&\hspace{0.5cm}
+
\bigg(-\frac{{\sf e}^{2}}{{\sf c}^{2}}
+
\frac{\overline{\mathcal{L}}(\bar{\kaux})}{\mathcal{L}_{1}(\bar{\kaux})}
\frac{{\sf e}\bar{\sf e}}{\bar{\sf c}^{2}}
+
\frac{\overline{\mathcal{L}}_{1}(\overline{\mathcal{L}}_{1}(\kaux))}{\overline{\mathcal{L}}_{1}(\kaux)}
\frac{\sf e}{{\sf c}\bar{\sf c}}\bigg)\kappa\wedge\overline{\kappa}
+
\bigg(
\frac{\sf e}{\sf c}
-
\frac{\overline{\mathcal{L}}_{1}(\bar{\kaux})}{\mathcal{L}_{1}(\bar{\kaux})}\frac{\bar{\sf e}{\sf c}}{\bar{\sf c}^{2}}
-
\frac{\overline{\mathcal{L}}_{1}(\overline{\mathcal{L}}_{1}(\kaux))}{\overline{\mathcal{L}}_{1}(\kaux)}\frac{1}{\bar{\sf c}}\bigg)\ \zeta\wedge\bar{\kappa}\\
&\hspace{0.5cm}
-\frac{\overline{\mathcal{L}}_{1}(\bar{\kaux})}{\mathcal{L}_{1}(\bar{\kaux})}\frac{\sf e}{\bar{\sf c}}\ \kappa\wedge\bar{\zeta}
+
\frac{{\sf c}\overline{\mathcal{L}}_{1}(\bar{\kaux})}{\bar{\sf c}\mathcal{L}_{1}(\bar{\kaux})}\ \zeta\wedge\bar{\zeta}.
\end{aligned}
\end{equation*}
As before, we proceed with the absorption by setting 
\begin{equation*}
\begin{aligned}
\alpha &= \hat{\alpha}-
x_{\rho}\rho - x_{\kappa}\kappa - x_{\zeta}\zeta - x_{\bar{\kappa}}\bar{\kappa}- x_{\bar{\zeta}}\bar{\zeta},\\
\delta &= \hat{\delta}
-y_{\rho}\rho - y_{\kappa}\kappa - y_{\zeta}\zeta - y_{\bar{\kappa}}\bar{\kappa}- y_{\bar{\zeta}}\bar{\zeta}.
\end{aligned}
\end{equation*}
The equations that need attention are 
\begin{equation*}
\begin{aligned}
x_{\bar{\kappa}}+\overline{x_{\kappa}}
&= 
-\frac{\bar{\Paux}}{\bar{\sf c}}
-\frac{\overline{\mathcal{L}}_{1}(\bar{\kaux})}{\mathcal{L}_{1}(\bar{\kaux})}
\frac{\bar{\sf e}{\sf c}}{\bar{\sf c}^{2}},\\
x_{\bar{\kappa}} &= \frac{\sf e}{\sf c},\\
x_{\bar{\kappa}}-\overline{x_{\kappa}}
&=
-\frac{\sf e}{\sf c}
+
\frac{\overline{\mathcal{L}}_{1}(\bar{\kaux})}{\mathcal{L}_{1}(\bar{\kaux})}\frac{\bar{\sf e}{\sf c}}{\bar{\sf c}^{2}}
+
\frac{\overline{\mathcal{L}}_{1}(\overline{\mathcal{L}}_{1}(\kaux))}{\overline{\mathcal{L}}_{1}(\kaux)}\frac{1}{\bar{\sf c}}.
\end{aligned}
\end{equation*}
For the linear equations to have solutions, one therefore has to make the following choice for ${\sf e}$:
\[
{\sf e} = 
\frac{\sf c}{\bar{\sf c}}\bigg(
-\frac{1}{3}\bar{\Paux}+\frac{1}{3}\frac{\overline{\mathcal{L}}_{1}(\overline{\mathcal{L}}_{1}(\kaux))}{\overline{\mathcal{L}}_{1}(\kaux)}\bigg).
\]

We remark as well that in \cite{Foo-Merker-2019}, a similar normalisation is done during second loop of the Cartan process where the following choice for ${\sf b}$ is made: 
\[
{\sf b}=-\isqrt\bar{\sf c}{\sf{e}}+
\frac{\isqrt}{3}
{\sf c}\bigg(
\frac{\overline{\mathcal{L}}_{1}\overline{\mathcal{L}}_{1}(\kaux)}{\overline{\mathcal{L}}_{1}(\kaux)}
-\overline{\Paux}\bigg),
\]
 so that when ${\sf b}=0$ due to rigidity assumption, the same expression for ${\sf e}$ is also obtained.  At this stage, we set 
\[
\Baux := 
-\frac{1}{3}\bar{\Paux}+\frac{1}{3}\frac{\overline{\mathcal{L}}_{1}(\overline{\mathcal{L}}_{1}(\kaux))}{\overline{\mathcal{L}}_{1}(\kaux)}.
\]

\section{Final loop}
We make another change of base coframe by setting 
\[
\zeta_{0}'= \hat{\zeta}_{0}+\Baux\kappa_{0}.
\]
The new transformation group becomes 
\[
\left(
\begin{matrix}
\rho\\
\kappa\\
\zeta
\end{matrix}
\right)
=
\left(
\begin{matrix}
{\sf c}\bar{\sf c} & 0 & 0\\
0 & {\sf c} & 0\\
0 & 0 & \frac{\sf c}{\bar{\sf c}}
\end{matrix}
\right)
\left(
\begin{matrix}
\rho_{0}\\
\kappa_{0}\\
\zeta_{0}'
\end{matrix}
\right),
\]
with the new Darboux-Cartan structure:
\begin{equation}\label{finloop-iniCD}
\begin{aligned}
d\rho_{0} &= 
\bigg(\Paux
+
\Baux\frac{\mathcal{L}_{1}(\kaux)}{\overline{\mathcal{L}}_{1}(\kaux)}\bigg)\ \rho_{0}\wedge\kappa_{0}
-
\frac{\mathcal{L}_{1}(\kaux)}{\overline{\mathcal{L}}_{1}(\kaux)}\ \rho_{0}\wedge\zeta_{0}'\\
&\hspace{0.5cm}+
\bigg(\overline{\Paux}+
\overline{\Baux}\frac{\overline{\mathcal{L}}_{1}(\bar{\kaux})}{\mathcal{L}_{1}(\bar{\kaux})}\bigg)\ \rho_{0}\wedge\bar{\kappa}_{0}
-
\frac{\overline{\mathcal{L}}_{1}(\bar{\kaux})}{\mathcal{L}_{1}(\bar{\kaux})}\ \rho_{0}\wedge\bar{\zeta}_{0}'\\
&\hspace{0.5cm}
+\isqrt \kappa_{0}\wedge\overline{\kappa}_{0}\\
&=
\bigg(\Paux-\frac{\overline{\Paux}\mathcal{L}_{1}(\kaux)}{3\overline{\mathcal{L}}_{1}(\kaux)}
+
\frac{\overline{\mathcal{L}}_{1}(\overline{\mathcal{L}}_{1}(\kaux))\mathcal{L}_{1}(\kaux)}{3\overline{\mathcal{L}}_{1}(\kaux)^{2}}\bigg)\ \rho_{0}\wedge\kappa_{0}
-\frac{\mathcal{L}_{1}(\kaux)}{\overline{\mathcal{L}}_{1}(\kaux)}\ \rho_{0}\wedge\zeta_{0}'\\
&\hspace{0.5cm}
+\bigg(\overline{\Paux}
-\frac{\Paux\overline{\mathcal{L}}_{1}(\overline{\kaux})}{3\mathcal{L}_{1}(\overline{\kaux})}
+
\frac{\mathcal{L}_{1}(\mathcal{L}_{1}(\overline{\kaux}))\overline{\mathcal{L}}_{1}(\overline{\kaux})}{3\overline{\mathcal{L}}_{1}(\overline{\kaux})^{2}}\bigg)\ 
\rho_{0}\wedge\overline{\kappa}_{0}
-
\frac{\overline{\mathcal{L}}_{1}(\overline{\kaux})}{\mathcal{L}_{1}(\overline{\kaux})}\ 
\rho_{0}\wedge\overline{\zeta}_{0}'\\
&\hspace{0.5cm}
+\isqrt\kappa_{0}\wedge\overline{\kappa}_{0}
\\
&=:
\Raux_{1}\ \rho_{0}\wedge\kappa+\Raux_{2}\ \rho_{0}\wedge\zeta_{0}' + \overline{\Raux}_{1}\ \rho_{0}\wedge\bar{\kappa}_{0}+\overline{\Raux}_{2}\ \rho_{0}\wedge\bar{\zeta}_{0}'
+\isqrt\kappa_{0}\wedge\bar{\kappa}_{0},\\
d\kappa_{0} &= 
-\frac{\mathcal{L}_{1}(\kaux)}{\overline{\mathcal{L}}_{1}(\kaux)}\ \kappa_{0}\wedge\zeta_{0}'
-\Baux\ \kappa_{0}\wedge\bar{\kappa}_{0}
+\zeta_{0}'\wedge\bar{\kappa}_{0}\\
&=
-\frac{\mathcal{L}_{1}(\kaux)}{\overline{\mathcal{L}}_{1}(\kaux)}\ \kappa_{0}\wedge\zeta_{0}'
+
\bigg(
\frac{\overline{\Paux}}{3}
-
\frac{\overline{\mathcal{L}}_{1}(\overline{\mathcal{L}}_{1}(\kaux))}{3\overline{\mathcal{L}}_{1}(\kaux)}\bigg)\ \kappa_{0}\wedge\overline{\kappa}_{0}
+
\zeta_{0}'\wedge\overline{\kappa}_{0}\\
&=:
\Kaux_{5}\ \kappa_{0}\wedge\zeta_{0}'+
\Kaux_{6}\ \kappa_{0}\wedge\bar{\kappa}_{0}
+
\zeta_{0}'\wedge\bar{\kappa}_{0},
\end{aligned}
\end{equation}

The $2$-form $d\zeta_{0}'$ requires a bit of computation, as will be seen in the proof of the  following 
\begin{Proposition}
One has
\begin{equation*}
\begin{aligned}
d\zeta_{0}'
&=
\bigg(-\Baux \frac{\mathcal{L}_{1}(\kaux)}{\overline{\mathcal{L}}_{1}(\kaux)}
+
\frac{\mathcal{L}_{1}(\overline{\mathcal{L}}_{1}(\kaux))}{\overline{\mathcal{L}}_{1}(\kaux)}
-
\frac{\mathcal{K}(\Baux)}{\overline{\mathcal{L}}_{1}(\kaux)}\bigg)\ \kappa_{0}\wedge\zeta_{0}'
+
\bigg(
-\Baux^{2}+\Baux\frac{\overline{\mathcal{L}}_{1}(\overline{\mathcal{L}}_{1}(\kaux))}{\overline{\mathcal{L}}_{1}(\kaux)}-\overline{\mathcal{L}}_{1}(\Baux)\bigg)\ \kappa_{0}\wedge\bar{\kappa}_{0}\\
&\hspace{0.5cm} 
+\bigg(\Baux-\frac{\overline{\mathcal{L}}_{1}(\overline{\mathcal{L}}_{1}(\kaux))}{\overline{\mathcal{L}}_{1}(\kaux)}-\overline{\Baux}\frac{\overline{\mathcal{L}}_{1}(\bar{\kaux})}{\mathcal{L}_{1}(\bar{\kaux})}\bigg)\ \zeta_{0}'\wedge\bar{\kappa}_{0}
+
\frac{\overline{\mathcal{L}}_{1}(\bar{\kaux})}{\mathcal{L}_{1}(\bar{\kaux})}\ \zeta_{0}'\wedge\bar{\zeta}_{0}'\\
&=:
\Zaux_{5}\ \kappa_{0}\wedge\zeta_{0}'
+
\Zaux_{6}\ \kappa_{0}\wedge\bar{\kappa}_{0}
+
\Zaux_{8}\ \zeta_{0}'\wedge\bar{\kappa}_{0}
+
\Zaux_{9}\ \zeta_{0}'\wedge\bar{\zeta}_{0}'.
\end{aligned}
\end{equation*}
\end{Proposition}
\begin{proof}
Using the transformation $\zeta_{0}'=\hat{\zeta}_{0}+\Baux\kappa_{0}$, the $2$-forms $d\hat{\zeta}_{0}$ and $d\kappa_{0}$ are expressed in terms of the new coframe $(\rho,\kappa_{0},\zeta_{0}')$ as 
\begin{equation*}
\begin{aligned}
d\hat{\zeta}_{0} &=\frac{\mathcal{L}_{1}(\overline{\mathcal{L}}_{1}(\kaux))}{\overline{\mathcal{L}}_{1}(\kaux)}\ \kappa_{0}\wedge\zeta_{0}'
-
\Baux\frac{\overline{\mathcal{L}}_{1}(\overline{\kaux})}{\mathcal{L}_{1}(\overline{\kaux})}\ \kappa_{0}\wedge\overline{\zeta}_{0}'
-
\bigg(
\frac{\overline{\mathcal{L}}_{1}(\overline{\mathcal{L}}_{1}(\kaux))}{\overline{\mathcal{L}}_{1}(\kaux)}+\overline{\Baux}\frac{\overline{\mathcal{L}}_{1}(\overline{\kaux})}{\mathcal{L}_{1}(\overline{\kaux})}\bigg)\zeta_{0}'\wedge\overline{\kappa}_{0}\\
&\hspace{0.5cm}
+
\bigg(\Baux\frac{\overline{\mathcal{L}}_{1}(\overline{\mathcal{L}}_{1}(\kaux))}{\overline{\mathcal{L}}_{1}(\kaux)}+
\Baux\overline{\Baux}\frac{\overline{\mathcal{L}}_{1}(\overline{\kaux})}{\mathcal{L}_{1}(\overline{\kaux})}\bigg)\ \kappa_{0}\wedge\overline{\kappa}_{0}
+
\frac{\overline{\mathcal{L}}_{1}(\overline{\kaux})}{\mathcal{L}_{1}(\overline{\kaux})}\ \zeta_{0}'\wedge\overline{\zeta}_{0}',
\end{aligned}
\end{equation*}
as well as 
\begin{equation*}
\begin{aligned}
d\kappa_{0} &= 
-\frac{\mathcal{L}_{1}(\kaux)}{\overline{\mathcal{L}}_{1}(\kaux)}\ \kappa_{0}\wedge\zeta_{0}'
-\Baux\ \kappa_{0}\wedge\overline{\kappa}_{0}
+
\zeta_{0}'\wedge\overline{\kappa}_{0}.
\end{aligned}
\end{equation*}
Moreover one has for the $1$-form $d\Baux$ the following expansion
\begin{equation*}
\begin{aligned}
d\Baux
=
\mathcal{T}(\Baux)\ \rho_{0}
+
\mathcal{L}_{1}(\Baux)\ \kappa_{0}
+
\mathcal{K}(\Baux)\ \zeta_{0}
+
\overline{\mathcal{L}}_{1}(\Baux)\ \overline{\kappa}_{0}
+
\overline{\mathcal{K}}(\Baux)\overline{\zeta}_{0}.
\end{aligned}
\end{equation*}
By rigidity assumption, $\mathcal{T}(\Baux)\equiv 0$; and by using the Assertion 7.4 on page 26 of Foo-Merker \cite{Foo-Merker-2019},
\[
\overline{\mathcal{K}}(\Baux)
=
-\Baux\overline{\mathcal{L}}_{1}(\overline{\kaux}).
\]
Using these two observations, the $1$-form $d\Baux$ is therefore
\begin{equation*}
\begin{aligned}
d\Baux &=
\bigg(\mathcal{L}_{1}(\Baux)-\Baux\frac{\mathcal{K}(\Baux)}{\overline{\mathcal{L}}_{1}(\kaux)}\bigg)\ \kappa_{0}
+
\frac{\mathcal{K}(\Baux)}{\overline{\mathcal{L}}_{1}(\kaux)}\ \zeta_{0}'
+
\bigg(\overline{\mathcal{L}}_{1}(\Baux)
+
\Baux\overline{\Baux}
\frac{\overline{\mathcal{L}}_{1}(\overline{\kaux})}{\mathcal{L}_{1}(\overline{\kaux})}\bigg)
\overline{\kappa}_{0}
-
\Baux
\frac{\overline{\mathcal{L}}_{1}(\overline{\kaux})}{\mathcal{L}_{1}(\overline{\kaux})}\ \overline{\zeta}_{0}'.
\end{aligned}
\end{equation*}
Substituting $d\hat{\zeta}_{0}$, $d\kappa_{0}$ and $d\Baux$ in the following identity 
\[
d\zeta_{0}'=d\hat{\zeta}_{0}+d\Baux\wedge\kappa_{0}+\Baux\ d\kappa_{0}
\]
by the expressions computed above finishes the proof of the proposition.
\end{proof}

Explicitly, 
\begin{equation}\label{finloop-iniCD2}
\begin{aligned}
d\zeta_{0}'
&=
\bigg(
\frac{\overline{\Paux}\mathcal{L}_{1}(\kaux)}{3\overline{\mathcal{L}}_{1}(\kaux)}
-
\frac{\overline{\mathcal{L}}_{1}(\overline{\mathcal{L}}_{1}(\kaux))\mathcal{L}_{1}(\kaux)}{3\overline{\mathcal{L}}_{1}(\kaux)^{2}}
+
\frac{\mathcal{L}_{1}(\overline{\mathcal{L}}_{1}(\kaux))}{\overline{\mathcal{L}}_{1}(\kaux)}
+
\frac{\mathcal{K}(\overline{\Paux})}{3\overline{\mathcal{L}}_{1}(\kaux)}\\
&\hspace{1cm}
-
\frac{\mathcal{K}(\overline{\mathcal{L}}_{1}(\overline{\mathcal{L}}_{1}(\kaux)))}{3\overline{\mathcal{L}}_{1}(\kaux)^{2}}
+
\frac{\mathcal{K}(\overline{\mathcal{L}}_{1}(\kaux))\overline{\mathcal{L}}_{1}(\overline{\mathcal{L}}_{1}(\kaux))}{3\overline{\mathcal{L}}_{1}(\kaux)^{3}}\bigg)\ \kappa_{0}\wedge\zeta_{0}'\\
&\hspace{0.5cm}
+\bigg(\frac{-\overline{\Paux}^{2}}{9}
-
\frac{\overline{\Paux}\overline{\mathcal{L}}_{1}(\overline{\mathcal{L}}_{1}(\kaux))}{9\overline{\mathcal{L}}_{1}(\kaux)}
+
\frac{5\overline{\mathcal{L}}_{1}(\overline{\mathcal{L}}_{1}(\kaux))^{2}}{9\overline{\mathcal{L}}_{1}(\kaux)^{2}}\\
&\hspace{1.5cm}
-
\frac{\overline{\mathcal{L}}_{1}(\overline{\mathcal{L}}_{1}(\overline{\mathcal{L}}_{1}(\kaux)))}{3\overline{\mathcal{L}}_{1}(\kaux)}
+
\frac{\overline{\mathcal{L}}_{1}(\overline{\Paux})}{3}\bigg)\kappa_{0}\wedge\overline{\kappa}_{0}\\
&\hspace{0.5cm}
+
\bigg(\frac{-\overline{\Paux}}{3}
-\frac{2\overline{\mathcal{L}}_{1}(\overline{\mathcal{L}}_{1}(\kaux))}{3\overline{\mathcal{L}}_{1}(\kaux)}
+
\frac{\Paux\overline{\mathcal{L}}_{1}(\overline{\kaux})}{3\mathcal{L}_{1}(\overline{\kaux})}
-
\frac{\mathcal{L}_{1}(\mathcal{L}_{1}(\overline{\kaux}))\overline{\mathcal{L}}_{1}(\overline{\kaux})}{3\mathcal{L}_{1}(\overline{\kaux})^{2}}\bigg)\zeta_{0}'\wedge\overline{\kappa}_{0}\\
&\hspace{0.5cm}
+
\frac{\overline{\mathcal{L}}_{1}(\overline{\kaux})}{\mathcal{L}_{1}(\overline{\kaux})}\ \zeta_{0}'\wedge\overline{\zeta}_{0}'.
\end{aligned}
\end{equation}

After transformation, the new $2$-forms $d\rho$, $d\kappa$ and $d\zeta$ become 
\begin{equation*}
\begin{aligned}
d\rho &= (\alpha+\bar{\alpha})\wedge\rho
+\frac{1}{\sf c}\Raux_{1}\ \rho\wedge\kappa
+
\frac{\bar{\sf c}}{\sf c}\Raux_{2}\ \rho\wedge\zeta
+
\frac{1}{\bar{\sf c}}\overline{\Raux}_{1}\ \rho\wedge\bar{\kappa}
+
\frac{\sf c}{\bar{\sf c}}\ \overline{\Raux}_{2}\ \rho\wedge\bar{\zeta}
+
\isqrt\kappa\wedge\bar{\kappa},\\
d\kappa &= \alpha\wedge\kappa
+
\frac{\bar{\sf c}}{\sf c}\Kaux_{5}\ \kappa\wedge\zeta
+
\frac{1}{\bar{\sf c}}\Kaux_{6}\ \kappa\wedge\bar{\kappa}
+
\zeta\wedge\bar{\kappa},\\
d\zeta &= 
(\alpha-\bar{\alpha})\wedge\zeta
+
\frac{1}{\sf c}\Zaux_{5}\ \kappa\wedge\zeta
+
\frac{1}{\bar{\sf c}^{2}}\Zaux_{6}\ \kappa\wedge\bar{\kappa}
+
\frac{1}{\bar{\sf c}}\Zaux_{8}\ \zeta\wedge\bar{\kappa}
+
\frac{\sf c}{\bar{\sf c}}\Zaux_{9}\ \zeta\wedge\bar{\zeta}.
\end{aligned}
\end{equation*}
By setting the new Maurer-Cartan 1-form as 
\[
\alpha := \hat{\alpha}
-
x_{\rho}\rho-x_{\kappa}\kappa-x_{\zeta}\zeta-x_{\bar{\kappa}}\bar{\kappa}-x_{\bar{\zeta}}\bar{\zeta},
\]
with 
\begin{equation}\label{abs-fin}
x_{\rho}=0,\qquad x_{\kappa}=-\frac{1}{\sf c}\Raux_{1}+\frac{1}{\sf c}\bar{\Kaux}_{6},\qquad
x_{\bar{\kappa}}=\frac{1}{\bar{\sf c}}\Baux,\qquad
x_{\zeta}=\frac{\bar{\sf c}\mathcal{L}_{1}(\kaux)}{{\sf c}\overline{\mathcal{L}}_{1}(\kaux)},\qquad
x_{\bar{\zeta}}=0,
\end{equation}
the final absorbed equations become:
\begin{eqnarray}
d\rho &=& (\hat{\alpha}+\overline{\hat{\alpha}})\wedge\rho + \isqrt \kappa\wedge\bar{\kappa},\label{final-1}\\
d\kappa &=& \hat{\alpha}\wedge\kappa + \zeta\wedge\bar{\kappa},\label{final-2}\\
d\zeta &=& (\hat{\alpha}-\overline{\hat{\alpha}})\wedge\zeta
+
\frac{1}{\sf c}(\Zaux_{5}-\overline{\Zaux}_{8})\ \kappa\wedge\zeta
+
\frac{1}{\bar{\sf c}^{2}}\Zaux_{6}\ \kappa\wedge\bar{\kappa}\label{final-3}.
\end{eqnarray}

\section{The $\{e\}$-structure.}

This time, for ease of notation, we write
\[
{\sf S}_{5}=\frac{1}{\sf c}(\Zaux_{5}-\bar{\Zaux}_{8}):=\frac{1}{\sf c}\Iaux_{0},\qquad
{\sf S}_{6} = \frac{1}{\bar{\sf c}^{2}}\Zaux_{6}:=\frac{1}{\bar{\sf c}^{2}}\Vaux_{0}.\]

If we write
\[
\psi:= -{\sf S}_{5}\zeta-{\sf S}_{6}\bar{\kappa},
\]
equation \eqref{final-3} may be written otherwise as
\[
d\zeta=
(\hat{\alpha}-\overline{\hat{\alpha}})\wedge\zeta+\psi\wedge\kappa.
\]
Based on the model case in Section 4, one should obtain for $d\hat{\alpha}$ the following:
\[
d\hat{\alpha} = \zeta\wedge\bar{\zeta}+\cdots,
\]
where the remaining terms are $2$-forms that vanish in the model case. Taking exterior derivatives of both sides of equations \eqref{final-1},  \eqref{final-2} and \eqref{final-3}:
\begin{equation}\label{streqfinal-1}
\begin{aligned}
0 &= (d\hat{\alpha}+d\overline{\hat{\alpha}})\wedge\rho,\\
0 &= (d\hat{\alpha}-\zeta\wedge\bar{\zeta}+{\sf S}_{5}\zeta\wedge\bar{\kappa})\wedge\kappa\\
0 &=
(d\hat{\alpha}-d\bar{\hat{\alpha}})\wedge\zeta
-
(\hat{\alpha}-\bar{\hat{\alpha}})\wedge d\zeta
+
d\psi\wedge\kappa
-
\psi\wedge\alpha\wedge\kappa.
\end{aligned}
\end{equation}
In the second equation of \eqref{streqfinal-1}, Cartan's lemma provides a $1$-form $A$ with 
\[
d\hat{\alpha} = \zeta\wedge\bar{\zeta}
-
{\sf S}_{5}\zeta\wedge\bar{\kappa}
+
A\wedge\kappa.
\]

To study $A$, write it as a formal linear combination of the $1$-forms with unknown coefficients:
\[
A= A_{\rho}\rho+A_{\kappa}\kappa+A_{\zeta}\zeta+A_{\hat{\alpha}}\hat{\alpha}+
A_{\bar{\kappa}}\bar{\kappa}+A_{\bar{\zeta}}\bar{\zeta}+A_{\bar{\hat{\alpha}}}\bar{\hat{\alpha}}.
\]
From the first equation of \eqref{streqfinal-1}, one obtains 
\[
A_{\bar{\zeta}}=\overline{\sf S}_{5},\qquad
A_{\zeta}=0,\qquad
A_{\bar{\kappa}}\text{ is real},\qquad
A_{\hat{\alpha}}=A_{\bar{\hat{\alpha}}}=0,
\]
and so
\[
d\hat{\alpha}=\zeta\wedge\bar{\zeta}
-
{\sf S}_{5}\zeta\wedge\bar{\kappa}
+
A_{\rho}\rho\wedge\kappa
+
A_{\bar{\kappa}}\bar{\kappa}\wedge\kappa
+
\overline{\sf S}_{5}\bar{\zeta}\wedge\kappa.
\]

Using this expression of $d\hat{\alpha}$ in the third equation of \eqref{streqfinal-1}, the remaining coefficients of $A$ are therefore obtained:
\[
A_{\rho}=0,\qquad 
0=2A_{\bar{\kappa}}\ \bar{\kappa}\wedge\kappa\wedge\zeta\wedge\bar{\zeta}
+
\bar{\alpha}\wedge\psi\wedge\kappa\wedge\bar{\zeta}
+
d\psi\wedge\kappa\wedge\bar{\zeta}.
\]
We expand $d\psi$ so that 
\begin{equation}
\begin{aligned}
d\psi \wedge\kappa\wedge\bar{\zeta}
&=
(-d{\sf S}_{5}\wedge\zeta-{\sf S}_{5}d\zeta-d{\sf S}_{6}\wedge\bar{\kappa}-{\sf S}_{6}d\bar{\kappa})\wedge\kappa\wedge\bar{\zeta}\\
&=
-(({\sf S}_{5})_{\bar{\kappa}}-
({\sf S}_{6})_{\zeta})\ \kappa\wedge\bar{\kappa}\wedge\zeta\wedge\bar{\zeta}+\cdots,
\end{aligned}
\end{equation}
where $(\bullet)_{\bar{\kappa}}$ denotes the covariant derivative of the function $\bullet$ with respect to $\bar{\kappa}$ (and same definition applies to $(\bullet)_{\zeta}$). We could have concluded the $\{e\}$-structure by declaring 
\[
A_{\bar{\kappa}}= -{\textstyle{\frac{1}{2}}}(({\sf S}_{5})_{\bar{\kappa}}-
({\sf S}_{6})_{\zeta}),
\]
which is a secondary invariant. 

To make sure that the equation does make sense, the term on the right needs to be verified that it is real-valued. This requires some computation. First we need a lemma:
\begin{Lemma}\label{lemma-10.3}
On the $G$-structure $M\times G^{2}$ with coordinates $(z_{1},z_{2},\overline{z}_{1},\overline{z}_{2},v,{\sf c},\overline{\sf c})$, let $F:M\times G^{2}\rightarrow \mathbb{C}$ be a function. Then 
\begin{equation}
\begin{aligned}
dF &= 
c\partial_{\sf c}F\ \hat{\alpha}
+
\bar{\sf c}\partial_{\bar{\sf c}}F\ \bar{\hat{\alpha}}
+
\bigg(
\frac{1}{{\sf c}\bar{\sf c}}\mathcal{T}(F)-
{\sf c}x_{\rho}\partial_{\sf c}F
-\bar{\sf c}\overline{x_{\rho}}\partial_{\bar{\sf c}}F\bigg)\rho\\
&\hspace{0.5cm}
+\bigg(
\frac{1}{\sf c}
\bigg(\mathcal{L}_{1}(F)-\Baux
\frac{\mathcal{K}(F)}{\overline{\mathcal{L}}_{1}(\kaux)}\bigg)
-
{\sf c}x_{\kappa}\partial_{\sf c}F
-\bar{\sf c}\overline{x_{\bar{\kappa}}}\partial_{\bar{\sf c}}F\bigg)\kappa\\
&\hspace{0.5cm}
+\bigg(
\frac{\bar{\sf c}}{\sf c}\frac{\mathcal{K}(F)}{\overline{\mathcal{L}}_{1}(\kaux)}
-{\sf c}x_{\zeta}\partial_{\sf c}F
-\bar{c}\overline{x_{\bar{\zeta}}}\partial_{\bar{c}}F\bigg)\zeta\\
&\hspace{0.5cm}
+\bigg(
\frac{1}{\bar{\sf c}}\bigg(\overline{\mathcal{L}}_{1}(F)-\overline{\Baux}
\frac{\overline{\mathcal{K}}(F)}{\mathcal{L}_{1}(\overline{\kaux})}\bigg)
-{\sf c}x_{\bar{\kappa}}\partial_{\sf c}F
-\bar{\sf c}\overline{x_{\kappa}}\partial_{\bar{\sf c}}F\bigg)\bar{\kappa}\\
&\hspace{0.5cm}+
\bigg(\frac{\sf c}{\bar{\sf c}}
\frac{\bar{\mathcal{K}}(F)}{\mathcal{L}_{1}(\bar{\kaux})}
-{\sf c}x_{\bar{\zeta}}\partial_{\sf c}F
-
\bar{\sf c}\overline{x_{\zeta}}\partial_{\bar{\sf c}}F\bigg)\bar{\zeta}\\
&:=
\partial_{\alpha}(F)\ \alpha
+
\partial_{\overline{\alpha}}(F)\ \overline{\alpha}
+
\partial_{\rho}(F)\ \rho
+
\partial_{\kappa}(F)\ \kappa
+
\partial_{\zeta}(F)\ \zeta\\
&\hspace{0.5cm}
+
\partial_{\overline{\kappa}}(F)\ \overline{\kappa}
+
\partial_{\overline{\zeta}}(F)\ \overline{\zeta}.
\end{aligned}
\end{equation}
\end{Lemma} 

The proof of the lemma is done by straightforward computation which will be skipped.  With the solution to the absorption equations \eqref{abs-fin}, we therefore have the following vector fields:

\begin{equation}\label{vect-abs-fin}
\begin{aligned}
\partial_{\alpha} &:= {\sf c}\partial_{\sf c},\\
\partial_{\rho} &:= \frac{1}{\sf c\overline{c}}\mathcal{T},\\
\partial_{\kappa} &:= 
\frac{1}{\sf c}\bigg(\mathcal{L}_{1}-\Baux\frac{\mathcal{K}}{\overline{\mathcal{L}}_{1}(\kaux)}\bigg)
-{\sf c}\bigg(-\frac{1}{\sf c}\Raux_{1}+\frac{1}{\sf c}\overline{\Kaux}_{6}\bigg)\partial_{\sf c}
+\frac{\overline{\sf c}}{\sf c}\overline{\Baux}\partial_{\overline{\sf c}},\\
\partial_{\zeta} &= \frac{\overline{\sf c}}{\sf c}\frac{\mathcal{K}}{\overline{\mathcal{L}}_{1}(\kaux)}
-
\overline{\sf c}\frac{\mathcal{L}_{1}(\kaux)}{\overline{\mathcal{L}}_{1}(\kaux)}\partial_{\sf c},
\end{aligned}
\end{equation}
while the vector fields $\partial_{\overline{\alpha}}$, $\partial_{\overline{\kappa}}$, $\partial_{\overline{\zeta}}$ are respective complex conjugates of $\partial_{\alpha}$, $\partial_{\kappa}$, $\partial_{\zeta}$. As a result:
\begin{equation*}
\begin{aligned}
({\sf S}_{5})_{\bar{\kappa}}
-
({\sf S}_{6})_{\zeta}
&=
\frac{1}{{\sf c}\bar{\sf c}}
\bigg(
\overline{\mathcal{L}}_{1}\big(\Iaux_{0}\big)
-\overline{\Baux}
\frac{\overline{\mathcal{K}}\big(\Iaux_{0}\big)}{\mathcal{L}(\overline{\kaux})}
+\Baux\Iaux_{0}
-
\frac{\mathcal{K}\big(\Vaux_{0}\big)}{\overline{\mathcal{L}}(\kaux)}
\bigg):=
\frac{1}{\sf c\overline{c}}\Qaux_{0}\\
&=
\frac{1}{{\sf c}\bar{\sf c}}
\bigg(
\overline{\mathcal{L}}_{1}(\Zaux_{5})
-
\overline{\mathcal{L}}_{1}(\overline{\Zaux}_{8})
-
\overline{\Baux}\frac{\overline{\mathcal{K}}(\Zaux_{5})}{\mathcal{L}_{1}(\overline{\kaux})}
+
\overline{\Baux}\frac{\overline{\mathcal{K}}(\overline{\Zaux}_{8})}{\mathcal{L}_{1}(\overline{\kaux})}
+
\Baux\Zaux_{5}
-
\Baux\overline{\Zaux}_{8}
-
\frac{\mathcal{K}(\Zaux_{6})}{\overline{\mathcal{L}}_{1}(\kaux)}\bigg).
\end{aligned}
\end{equation*}
We will also need the following
\begin{Lemma}
One has the following identity
\begin{equation*}
\begin{aligned}
\overline{\mathcal{L}}_{1}(\Zaux_{5})-\frac{\mathcal{K}(\Zaux_{6})}{\overline{\mathcal{L}}_{1}(\kaux)}
&=
\overline{\Baux}
\frac{\overline{\mathcal{K}}(\Zaux_{5})}{\mathcal{L}_{1}(\bar{\kaux})}
+
\Zaux_{5}\Kaux_{6}
-\Zaux_{6}\Kaux_{5}
-
\mathcal{L}_{1}(\Zaux_{8})
+\Baux\frac{\mathcal{K}(\Zaux_{8})}{\overline{\mathcal{L}}_{1}(\kaux)}
+
\Zaux_{8}\overline{\Kaux}_{6}
+
\Zaux_{9}\overline{\Zaux}_{6}.
\end{aligned}
\end{equation*}
\end{Lemma}
\begin{proof}
We will compute the terms on the left-hand side by applying $d^{2}\equiv 0$ to the third equation of equation \eqref{finloop-iniCD}. Doing so, while wedging on both sides of $d^{2}\zeta_{0}'=0$ with $\rho\wedge\bar{\zeta}_{0}'$, one should get
\begin{equation*}
\begin{aligned}
0
=
\big(
(\Zaux_{5})_{\bar{\kappa}_{0}}
-
\Zaux_{5}\Kaux_{6}
-
\Zaux_{5}\Zaux_{8}
-
(\Zaux_{6})_{\zeta_{0}'}
+
\Zaux_{6}\Kaux_{5}
+
(\Zaux_{8})_{\kappa_{0}}
+
\Zaux_{8}\Zaux_{5}
-
\Zaux_{8}\overline{\Kaux}_{6}
-
\Zaux_{9}\overline{\Zaux}_{6}\big)
\rho_{0}\wedge\kappa_{0}\wedge\bar{\kappa}_{0}\wedge\zeta'_{0}\wedge\bar{\zeta}'_{0}.
\end{aligned}
\end{equation*}
Finally, for any function $G$ independent of ${\sf c}$, one uses the following formula
\begin{equation*}
\begin{aligned}
dG &= \mathcal{T}(G)\rho
+
\bigg(\mathcal{L}_{1}(G)-\Baux\frac{\mathcal{K}(G)}{\overline{\mathcal{L}}_{1}(\kaux)}\bigg)
\kappa_{0}
+
\frac{\mathcal{K}(G)}{\overline{\mathcal{L}}_{1}(\kaux)}\zeta_{0}'\\
&\hspace{0.5cm}
+\bigg(\overline{\mathcal{L}}_{1}(G)
-
\overline{\Baux}
\frac{\overline{\mathcal{K}}(G)}{\overline{\mathcal{L}}_{1}(\kaux)}\bigg)
\overline{\kappa}_{0}
+
\frac{\overline{\mathcal{K}}(G)}{\overline{\mathcal{L}}_{1}(\kaux)}\overline{\zeta}_{0}'.
\end{aligned}
\end{equation*}
The proof is therefore complete by applying this to $(\Zaux_{5})_{\bar{\kappa}_{0}}$, $(\Zaux_{6})_{\zeta_{0}'}$ and $(\Zaux_{8})_{\kappa_{0}}$.
\end{proof}
Substituting the identity into $A_{\bar{k}}$, one has therefore
\begin{equation}\label{Q0isreal}
\begin{aligned}
-2A_{\bar{\kappa}}
&=
\frac{1}{{\sf c}\bar{\sf c}}
\bigg(
(-\Zaux_{6}\Kaux_{5}+\Zaux_{9}\overline{\Zaux}_{6})
-
\mathcal{L}_{1}(\Zaux_{8})
-
\overline{\mathcal{L}}_{1}(\overline{\Zaux_{8}})
+
\Baux\frac{\mathcal{K}(\Zaux_{8})}{\overline{\mathcal{L}}_{1}(\kaux)}
+
\overline{\Baux}\frac{\overline{\mathcal{K}}(\overline{\Zaux}_{8})}{\mathcal{L}_{1}(\overline{\kaux})}
-
\Zaux_{8}\overline{\Baux}
-
\overline{\Zaux}_{8}\Baux
\bigg),
\end{aligned}
\end{equation}
and observing that $\Zaux_{9}=-\overline{\Kaux}_{5}$, the coefficient $A_{\bar{\kappa}}$ is thus real-valued, and  the $\{e\}$-structure is finally complete. 

We have therefore proved Theorem \ref{thm-main-equiv-2}.

In the interest of computations, the secondary invariant 
\[
\Qaux_{0}
:=
\frac{1}{2}
\bigg(
\overline{\mathcal{L}}_{1}\big(\Iaux_{0}\big)
-\overline{\Baux}
\frac{\overline{\mathcal{K}}\big(\Iaux_{0}\big)}{\mathcal{L}_{1}(\overline{\kaux})}
+\Baux\Iaux_{0}
-
\frac{\mathcal{K}\big(\Vaux_{0}\big)}{\overline{\mathcal{L}}_{1}(\kaux)}
\bigg)
\]
may further be simplified using the following:

\begin{Proposition}
Under the Levi degeneracy assumption, one has:
\[
\frac{\overline{\mathcal{K}}(\Iaux_{0})}{\mathcal{L}_{1}(\overline{\kaux})}
=
-2\overline{\Iaux_{0}}.\]
\end{Proposition}

\begin{proof}
We remark here that the Levi-degeneracy condition is necessary to normalise the expression and thus it cannot be dropped. It is implicitly used in $d\rho$ the first equation of the following $\{e\}$-structure:
\begin{equation*}
\begin{aligned}
d\rho &= (\alpha+\overline{\alpha})\wedge\rho + \isqrt \kappa\wedge\overline{\kappa},\\
d\kappa &= \alpha\wedge\kappa + \zeta\wedge\overline{\kappa},\\
d\zeta &= (\alpha-\overline{\alpha})\wedge\zeta
+
\frac{1}{\sf c}\Iaux_{0}\ \kappa\wedge\zeta
+
\frac{1}{\overline{\sf c}\overline{\sf c}}
\Vaux_{0}\ \kappa\wedge\overline{\kappa},\\
d\alpha &= 
\zeta\wedge\overline{\zeta}
-
\frac{1}{\sf c}\Iaux_{0}\ \zeta\wedge\overline{\kappa}
+
\frac{1}{{\sf c}\overline{\sf c}}\Qaux_{0}\ \kappa\wedge\overline{\kappa}
+
\frac{1}{\overline{\sf c}}\overline{\Iaux}_{0}\ \overline{\zeta}\wedge\kappa.
\end{aligned}
\end{equation*}
Applying Poincar\'{e} derivative to the third equation $d\zeta$ and using $d^2\equiv 0$, while wedging on both sides with $\alpha\wedge\overline{\alpha}\wedge\rho\wedge\overline{\kappa}$, we obtain 
\begin{equation*}
\begin{aligned}
0 &= d\alpha\wedge\zeta\wedge\alpha\wedge\overline{\alpha}\wedge\rho\wedge\overline{\kappa}
-d\overline{\alpha}\wedge\zeta\wedge\alpha\wedge\overline{\alpha}\wedge\rho\wedge\overline{\kappa}\\
&\hspace{0.5cm} +
\partial_{\overline{\zeta}}\bigg(\frac{1}{{\sf c}}\Iaux_{0}\bigg)\overline{\zeta}\wedge\kappa\wedge\zeta\wedge\alpha\wedge\overline{\alpha}\wedge\rho\wedge\overline{\kappa},
\end{aligned}
\end{equation*}
where $\partial_{\overline{\zeta}}$ is the following vector field coming from equation \eqref{vect-abs-fin}:
\[
\partial_{\overline{\zeta}}
=
\frac{\sf c}{\overline{\sf c}}\frac{\overline{\mathcal{K}}}{\mathcal{L}_{1}(\overline{\kaux})}
-
{\sf c}
\frac{\overline{\mathcal{L}}_{1}(\overline{\kaux})}{\mathcal{L}_{1}(\overline{\kaux})}\partial_{\overline{\sf c}}.
\]
Then using $d\alpha$ and $d\overline{\alpha}$ from the $\{e\}$-structure, we obtain the desired identity.
\end{proof}

Thus we recover the expression of $\Qaux_{0}$ as appeared in the introduction.

\bigskip
{\scriptsize
{\sc Wei Guo {\sc Foo}. Hua Loo-Keng Center for Mathematical Sciences, Academy of Mathematics and Systems Science, Chinese Academy of Sciences, Beijing, China.\\
Email address:} \texttt{fooweiguo@hotmail.com}}\\[-10pt]

{\scriptsize
{\sc Jo\"el {\sc Merker}. Laboratoire de Math\'{e}matiques d'Orsay, Universit\'{e} Paris-Sud, CNRS, Universit\'{e} Paris-Saclay, 91405 Orsay Cedex, France.\\
Email address:} \texttt{joel.merker@math.u-psud.fr}}\\[-10pt]

{\scriptsize
{\sc The-Anh {\sc Ta}. Laboratoire de Math\'{e}matiques d'Orsay, Universit\'{e} Paris-Sud, CNRS, Universit\'{e} Paris-Saclay, 91405 Orsay Cedex, France.\\
Email address:} \texttt{the-anh.ta@u-psud.fr}}


\vfill\end{document}